\documentclass{article}
\usepackage[numbers]{natbib}
\usepackage[preprint]{neurips_2024}

\usepackage[utf8]{inputenc} 
\usepackage[T1]{fontenc}    
\usepackage{hyperref}       
\usepackage{url}            
\usepackage{booktabs}       
\usepackage{amsfonts}       
\usepackage{nicefrac}       
\usepackage{microtype}      
\usepackage{xcolor}         

\usepackage{amssymb}   
\usepackage{algorithm}
\usepackage{algpseudocode}  
\usepackage{amsthm}
\usepackage{graphicx}
\usepackage{stmaryrd}
\usepackage{nccmath}
\usepackage{geometry}
\usepackage{amsmath}
\newcommand{\R}{\mathbb{R}}
\newcommand{\C}{\mathbb{C}}
\newcommand{\N}{\mathbb{N}}
\newcommand{\E}{\mathbb{E}}
\newcommand{\V}{\mathbb{V}}

\newcommand{\tr}{\text{Tr}}
\newcommand{\rk}{\text{rank}}

\newcommand{\p}{\mathbb{P}}

\usepackage{physics}

\DeclareMathOperator*{\Diag}{Diag}

\newtheorem{thm}{Theorem}

\newtheorem{crl}{Corollary}
\newtheorem{dft}{Definition}
\newtheorem{lem}{Lemma}
\newtheorem{as}{Assumption}
\newtheorem{notat}{Notation}
\newtheorem{rmk}{Remark}
\newtheorem{ex}{Example}

\title{Asymptotic spectrum of weighted sample covariance: another proof of spectrum convergence}

\author{
  Benoît Oriol \\
  CEREMADE, Université Paris-Dauphine, PSL, Paris, France \\
Société Générale Corporate and Investment Banking, Puteaux, France \\
  \texttt{benoit.oriol@dauphine.eu} \\
}

\begin{document}

\maketitle

\begin{abstract}
We propose another proof of the high dimensional spectrum convergence of the weighted sample covariance, more concise and self-sufficient but with stronger, but reasonable assumptions. We explain and illustrates this theorem for different weight distributions and show how the spectrum behaves in finite samples with heavy tails. The general purpose is to provide a detailed introduction to the high dimensional spectrum of weighted sample covariance.
\end{abstract}

\phantomsection{Keywords:} weighted covariance, asymptotic spectrum, Random Matrix Theory, Cauchy-Stieltjes transform

\section{Introduction and related work}
The spectrum of sample covariance and its behavior in Kolmogorov asymptotics, when the dimension grows approximately linearly with the number of samples, is a deeply studied topic in high dimensional covariance and precision matrix estimation. Among the large literature on covariance estimation, the high dimensional framework has known several recent improvements thanks to Random Matrix Theory (RMT). Main high dimensional approaches are gathered in recent surveys of the field \cite{Pourahmadi2011, Fan2016, Ke2019, Ledoit2020a, Clifford2020}.

A central result on the sample covariance in high dimension is the Marcenko-Pastur theorem \cite{Marcenko1967}: it states that in high dimension, the spectrum of the sample covariance converges to a non-random distribution, characterized by an equation on its Cauchy-Stieltjes transfrom. This can be seen as a law of large number for the sample spectrum. Several works followed to weaken the assumptions, make the proof simpler, and sutdy the behavior of the asymptotic distribution \cite{Silverstein1995a, Silverstein1995b, Silverstein1995c, Bai1998}. Later, a central limit theorem for linear spectral statistics was also proven \cite{Bai2004}, and recently used for a statistical test on covariance structure \cite{Bodnar2024}.

For applications such as wireless communications for MIMO systems \cite{Couillet2011}, neurosciences for dynamic brain connectivity \cite{Honnorat2022}, finance with stochastic variance processes \cite{Bun2016, Bongiorno2023}, the \textbf{weighted} sample covariance naturally emerges. While the sample covariance is of the form $S = \sqrt{\Sigma} X X^* \sqrt{\Sigma}$ with $\Sigma$ the population covariance and $X$ the normalized noise matrix, the weighted sample covariance has a more general form $S = \sqrt{\Sigma} X W X^* \sqrt{\Sigma}$ with $W$ a weight matrix, typically nonnegative semi-definite.

The Marcenko-Pastur theorem was generalized to this setting of weighted sample covariance: in high dimension, the spectrum of a weighted sample covariance converges almost surely to a distribution solving a Marcenko-Pastur-like functional equation on its Cauchy-Stieltjes transform. This result appears under different sets of assumptions in different works \cite{Monvel1996, Tulino2004, Burda2005,Zhang2007, Paul2009}, the weaker set of assumptions being found in \cite{Zhang2007}.

This work has the pedagogical purpose of explaining and illustrating the theorem and providing a concise proof, at least compared to \cite{Zhang2007}, and mostly self-sufficient of the convergence of the spectrum of weighted sample covariance, for a set assumption stronger than the one found in \cite{Zhang2007} but still arguably quite general.

\section{Notation and assumptions}
Let us introduce the following notation. Notation is not constant across major works on the spectrum of sample covariances. In our work, we chose to follow mostly Silverstein one \cite{Silverstein1995a}.

\begin{notat}[The data matrix]
	There are $N$ samples of dimension $n$. We have:
	 \begin{itemize}
	 	\item $c_n = \frac{n}{N}$ the concentration ratio,
	 	\item $Z_n$ is the noise $n \times N$ matrix composed of i.i.d centered complex entries of variance $1$, 
	 	\item $T_n$ is the true covariance, a non-negative definite Hermitian matrix of size $n \times n$, 
	 	\item $W_n$ is the weight matrix, a $N \times N$ diagonal non-negative real matrix, 
		\item $Y_n = T_n^{1/2} Z_n $ is the observed data matrix.  
	\end{itemize}
	In the following, the subscripts $n$ are omitted when no confusion is possible.
\end{notat}

In this work, the object of interest is the weighted sample covariance $B_n$, particularly its spectrum of eigenvalues which are introduced below.
\begin{notat}[Weighted sample covariance]
	For $n \in \N^*$, the weighted sample covariance is defined by:
	\begin{equation}\label{}
	\begin{aligned}
		B_n := \frac{1}{N} Y_n W_n Y_n^*.
	\end{aligned}
	\end{equation} 
	We note $(\tau_1^{(n)},...,\tau_n^{(n)})$ the eigenvalues of $T_n$ in decreasing order, and $(\lambda_1^{(n)},...,\lambda_n^{(n)})$ the eigenvalues of $B_n$ in decreasing order.
	
	Additionally, we will need matrix in the proof:
	\begin{equation}\label{}
	\begin{aligned}
		\underline{B}_n := \frac{1}{N} W_n^{1/2} Y_n^* Y_n W_n^{1/2}.
	\end{aligned}
	\end{equation}
	$\underline{B}_n$ and $B_n$ have the same eigenvalues, except for $|n-N|$ zero eigenvalues.
\end{notat}

\begin{ex}[Standard weighting]
	The most common choice of weighting is of course the constant weighting: $W_n = I_N$. In this situation, $B_n$ is the standard sample covariance, and its asymptotic spectrum is covered in RMT by the work of Marcenko and Pastur \cite{Marcenko1967}.
\end{ex}

\begin{ex}[Exponentially weighted scheme]
	Another common choice in time series analysis is the exponentially weighted scheme. Parametrized by some $\alpha \in \R_+^*$, we define the weights as:
	\begin{equation}\label{}
	\begin{aligned}
		\forall i \in \llbracket 1,N \rrbracket, (W_n)_{ii} = \beta e^{-\alpha i/N}, \\
		\beta = e^{-\alpha /N} \frac{1 - e^{-\alpha /N}}{1 - e^{-\alpha}}.
	\end{aligned}
	\end{equation}
\end{ex}

\begin{notat}[Empirical spectrum distribution]
	We consider a Hermitian matrix $A$ of size $n \times n$ with real eigenvalues $(\mu_1,...,\mu_n)$. We define the empirical spectrum distribution of $A$, $F^A$, as:
	\begin{equation}\label{}
	\begin{aligned}
		F^A := \frac{1}{n} \sum_{i=1}^n \mathrm{1}_{[\mu_i, +\infty[}.
	\end{aligned}
	\end{equation}
	In particular, we have:
	\begin{equation}\label{}
	\begin{aligned}
		&F_n := F^{B_n} = \frac{1}{n} \sum_{i=1}^n \mathrm{1}_{[\lambda_i^{(n)}, +\infty[}, \\
		&\underline{F}_n := F^{\underline{B}_n} = (1 - c_n)\mathrm{1}_{[0, +\infty[} + c_n F_n.
	\end{aligned}
	\end{equation}
\end{notat}

A key tool to manipulate and understand the asymptotic behavior of $F_n$ is the Cauchy-Stieltjes transform: a transform of a real finite measure into a complex function defined on  $\C_+ = \{ z \in \C | \Im(z) > 0 \}$.
\begin{notat}[Cauchy-Stieltjes transform]
	For a real finite measure $H$, we denote $m_H$ its Cauchy-Stieltjes transform: 
	\begin{equation}\label{}
	\begin{aligned}
		\forall z \in \C_+, m_H(z) = \int \frac{1}{\tau - z} dH(\tau).
	\end{aligned}
	\end{equation}
	In particular, for all $z \in \C_+$, we denote:
	\begin{equation}\label{}
	\begin{aligned}
		m_n(z) := m_{F_n}(z) &= \frac{1}{n} \sum_{i=1}^n \frac{1}{ \lambda_i^{(n)}-z} \\
		&= \frac{1}{n} \tr \left((B_n - zI)^{-1} \right), \\
		\underline{m}_n(z) := m_{\underline{F}_n}(z) &= \frac{1}{N} \tr \left((\underline{B}_n - zI)^{-1} \right) \\
		&= - \frac{1-c_n}{z} + c_n m_n(z).
	\end{aligned}
	\end{equation}
\end{notat}

We handle the weakly convergence of $F_n$ through the pointwise almost sure convergence of its Cauchy-Stieltjes transform $m_n$ using the following fundamental theorem.
\begin{thm}[Cauchy-Stieltjes convergence, Theorem 5.8.3 \cite{Fleermann2022}]
	Let $(\mu_n)_n$ be random probability measures, $\mu$ be a deterministic probability measure, and $Z \subset \C_+$ that has an accumulation point in $\C_+$, then with $\underset{n \rightarrow \infty}{\implies}$ denoting the weak convergence:
	\begin{equation}\label{}
	\begin{aligned}
		&\mu_n \underset{n \rightarrow \infty}{\implies} \mu \text{ almost surely} \iff \\
		&\forall z \in Z, m_{\mu_n}(z) \underset{n \rightarrow \infty}{\longrightarrow} m_\mu(z) \text{ almost surely}.
	\end{aligned}
	\end{equation}
	 
\end{thm}

We consider the following set of mild assumptions necessary to the proof of the main result. Except for Assumption (d) which concerns specifically the weights, those assumptions were introduced by Silverstein \cite{Silverstein1995c} in his proof of Marcenko-Pastur theorem for Hermitian matrices.
\begin{as}
	Assume that:
	\begin{itemize}
		\item[(a)] For all $n \in \N^*$ and $i \in \llbracket 1, N \rrbracket, j \in \llbracket 1, n \rrbracket$, the $(X_n)_{i,j} \in \C$ are i.i.d random variables with $\E[(X_n)_{i,j}] = 0$ and $\E[|(X_n)_{i,j}|^2] = 1$.
		\item[(b)] $N = N(n)$ and $c_n = \frac{n}{N(n)} \underset{n \rightarrow \infty}{\longrightarrow} c \in \R_+^*$.
		\item[(c)] $F^{T_n} \underset{n \rightarrow \infty}{\implies} H$ almost surely where $H$ is a probability distribution function - p.d.f. - on $\R_+$.
		\item[(d)] $F^{W_n} \underset{n \rightarrow \infty}{\implies} D$ almost surely where $D \in \mathcal{L}^1(\R_+)$ is a p.d.f., and almost surely $\int x dF^{W_n}(x) \underset{n \rightarrow \infty}{\longrightarrow} \int x dD(x) \in \R_+$.
		\item[(e)] $X_n$, $T_n$ and $W_n$ are independent.
	\end{itemize}
\end{as}
Assumption (a) makes sure that the covariance is well-defined and that $\E\left[\frac{1}{N} Y_n Y_n^*\right] = T_n$. Assumption (b) is the standard asymptotic setting we consider in RMT: the dimension and the number of samples are of the same order of magnitude. Assumption (c) and (d) ensure that the population covariance spectrum and the weight distribution converge: those two objects are essential to describe asymptotically the spectrum of $B_n$. The last assumption makes sure, for example, that the weights have no dependence on the drawn samples.

\section{Convergence of $m_n$ and $F_n$}
The convergence result egenralizing the Marcenko-Pastur theorem \cite{Marcenko1967} to weighted sample covariance is the following theorem. It characterizes the asymptotic non random spectrum $F$ of $B_n$ as $n \rightarrow \infty$ through its Cauchy-Stieltjes transform $m_F$. This result appears under different notation and sets of assumptions in different works \cite{Monvel1996, Tulino2004, Burda2005,Zhang2007, Paul2009}, the weaker set of assumptions being found in \cite{Zhang2007}. The set of assumptions (a)-(e) are not optimal, but we will be able to provide a concise proof, at least regarding \cite{Zhang2007}, still in a reasonably general setting.

\begin{thm}[Asymptotic spectrum]\label{dFE}
	Assume (a) to (e). Then, almost surely, $F_n$ converges weakly to a nonrandom p.d.f $F$, whose Cauchy-Stieltjes transform $m := m_F$ satisfies for all $z \in \C_+$:
	\begin{equation}\label{m}
	\begin{aligned}
		m(z) = \int \frac{1}{\tau \tilde m(z) - z}dH(\tau),
	\end{aligned}
	\end{equation} 
	where for all $z \in \C_+$, $\tilde m(z)$ is the unique solution in $\C \backslash \C_+$ of the following equation:
	\begin{equation}\label{}
	\begin{aligned}
		\tilde m(z) &=\int \frac{\delta}{1 + \delta c \int \frac{\tau}{\tau \tilde m(z) - z}dH(\tau)}dD(\delta).
	\end{aligned}
	\end{equation} 
\end{thm}
This theorem describes the asymptotic non-random behavior of $F_n$, through its Cauchy-Stieltjes transform, as Marcenko-Pastur theorem does for the standard sample covariance. 

\begin{rmk}[Role of $\tilde m$]
In Equation \eqref{m}, we note that if we fix $\tilde m(z) = 1$, then the right side of the equation becomes $m_H(z)$, so $F = H$. An interpretation of this equation is that $\tilde m(z)$ carries out the deformation of $H$ in its Cauchy-Stieltjes transform.
\end{rmk}

We only give here the main idea of the proof which is given in full detail in the Appendix. We follow the procedure introduced in the proof of \cite{Silverstein1995c}, and we make it work in our setting using recent developments of Lebesgue theorem \cite{Feinberg2019}. We firstly truncate the noise, the population eigenvalues and the weights by roughly $\log(n)$. This transform does not affect the asymptotic behavior of $F_n$. Doing that, we can apply a result from concentration theory that leads to the almost sure pointwise convergence of $m_n$ towards $m$. Through careful use of Lebesgue theorem for varying measures, we can extract an integral description of this limitnig $m$ that uses $\tilde m$. Finally, we check that $m$ is indeed the Cauchy-Stieltjes transform of some p.d.f., that we denote $F$.

\begin{figure}[]
\centering
\includegraphics[width=0.7\linewidth]{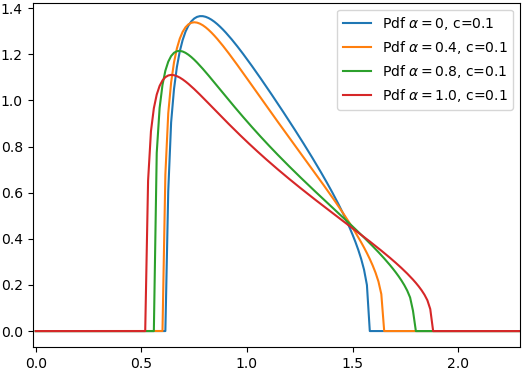}
\caption{Different approximations regarding parameter $\alpha$ of the asymptotic density of weighted sample eigenvalues for $c=0.1$, true covariance spectrum distribution $H = 1_{[1,\infty[}$ and weight distribution $D_\alpha = \frac{1}{2} 1_{[1-\alpha,\infty[} + \frac{1}{2} 1_{[1+\alpha,\infty[}$. $\alpha = 0$ corresponds to the classic Marcenko-Pastur density for $c=0.1$.}
\label{fig:sp_2diracs}
\end{figure}

\begin{figure}[]
\centering
\includegraphics[width=0.6\linewidth]{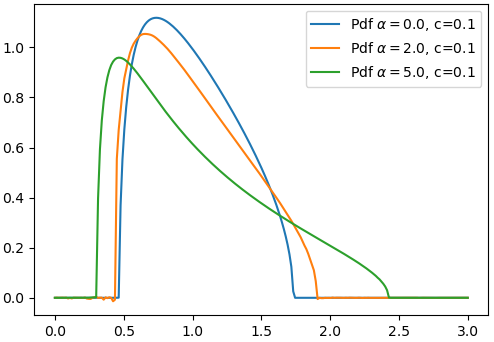}
\caption{Different approximations regarding parameter $\alpha$ of the asymptotic density of weighted sample eigenvalues for $c=0.1$, true covariance spectrum distribution $H = 1_{[1,\infty[}$ and weight distribution $D_\alpha: x \in [\beta e^{-\alpha}, \beta] \mapsto 1 + \frac{1}{\alpha}\log\left(\frac{x}{\beta}\right) \text{ with } \beta = \frac{\alpha}{1- e^{-\alpha}}$, defined in Definition \ref{ewma}. $\alpha = 0$ corresponds to the classic Marcenko-Pastur density for $c=0.1$.}
\label{fig:sp_exp}
\end{figure}

Examples of asymptotic densities are shown in figure \ref{fig:sp_2diracs} and \ref{fig:sp_exp}.

To understand better what is behind $\tilde m$, we can go back in the setting of the standard non-weighted covariance matrix. 

\begin{rmk}[Role of $\tilde m$ when $D =  1_{[1,+\infty[}$]
When $D =  1_{[1,+\infty[}$, we have that for all $z \in \C_+$, $\tilde m(z) = 1 - c\left(1+zm(z)\right) = -z \underline{m}(z)$.

The conclusion of this remark is that $\tilde{m}$ in the setting of weighted sample covariance has a similar role than $\underline{m}$ has in the uniformly weighted setting. 
\end{rmk}


\section{Experimental results}
\subsection{Impact of weight distribution on the asymptotic spectrum and on its support}
We show the behavior of the asymptotic density when we vary the weight distribution. The reference distribution of true eigenvalues is $H = 0.2 \times 1_{[1,\infty[} + 0.4 \times 1_{[3,\infty[} + 0.4 \times 1_{[10,\infty[}$, introduced by Bai and Silverstein \cite{Bai1998} and used by Ledoit and Péché \cite{Ledoit2009}. 

This distribution helps visualizing the phenomenon of spectral separation, which is deeply studied in \cite{Couillet2015}: areas of exclusion between the true eigenvalues where the asymptotic density is null. Those areas vanish when $c$ is increasing, as was shown in the classic setting of Marcenko-Pastur theory, but we can see here that the same phenomenon appears when the weight distribution is smoothly spreading.

For that matter, we are considering different weight distribution parametrized by $\alpha \in \R_+$ where $\alpha$ controls how much we are spreading the weights. We define the exponentially-weighted distribution, that corresponds to the distribution of weights in an Exponentially-Weighted Moving Average - EWMA -, used in time series analysis.

\begin{dft}[$\alpha$-exponentially weighted distribution]\label{ewma}
We fix $\alpha \in \R_+$. We define this law such as its cdf $D$ follows: $D(\beta e^{-\alpha t}) = 1 - t$ for $t \in [0,1]$. Moreover, we impose that $\int \delta dD(\delta) = 1$. We finally have: $\forall x \in [\beta e^{-\alpha}, \beta], D(x) = 1 + \frac{1}{\alpha} \log\left(\frac{x}{\beta}\right)$, with $\beta = \frac{\alpha}{1 - e^{-\alpha}}$.
\end{dft}

\begin{figure}[]
\centering
\includegraphics[width=0.7\linewidth]{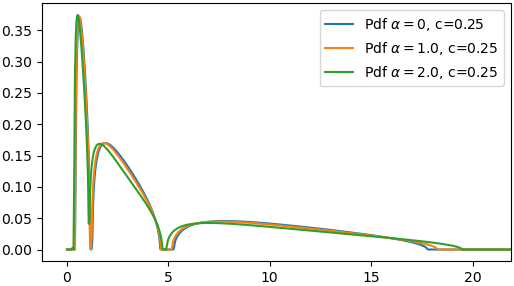}
\caption{Approximations of the asymptotic density of weighted sample eigenvalues for $c=0.25$,with $H = 0.2 \times 1_{[1,\infty[} + 0.4 \times 1_{[3,\infty[} + 0.4 \times 1_{[10,\infty[}$ and uniform weight distribution of parameter $\alpha$.}
\label{fig:sep_unif}
\end{figure}

\begin{figure}[]
\centering
\includegraphics[width=0.7\linewidth]{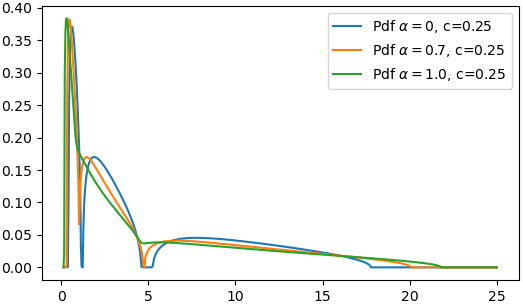}
\caption{Approximations of the asymptotic density of weighted sample eigenvalues for $c=0.25$, with $H = 0.2 \times 1_{[1,\infty[} + 0.4 \times 1_{[3,\infty[} + 0.4 \times 1_{[10,\infty[}$ and weight distribution from a mixture of 2 diracs of parameter $\alpha$.}
\label{fig:sep_2diracs}
\end{figure}

\begin{figure}[]
\centering
\includegraphics[width=0.7\linewidth]{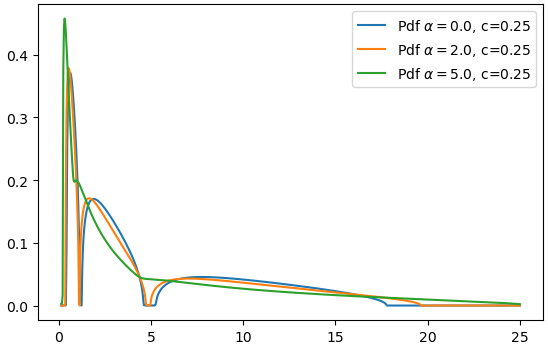}
\caption{Approximations of the asymptotic density of weighted sample eigenvalues for $c=0.25$, with $H = 0.2 \times 1_{[1,\infty[} + 0.4 \times 1_{[3,\infty[} + 0.4 \times 1_{[10,\infty[}$ and exponentially-weighted distribution of parameter $\alpha$.}
\label{fig:sep_exp}
\end{figure}

We consider three weight distributions:
\begin{itemize}
	\item A uniform distribution with $[1 - \alpha/2, 1 + \alpha/2]$ for $\alpha \in [0,2]$. The experiment is shown in figure \ref{fig:sep_unif}. There are 2 spectral separation for $\alpha = 0$ and $\alpha=1$. For $\alpha = 2$ there is only one left.
	\item A mixture of 2 diracs, with $D = \frac{1}{2} 1_{[1-\alpha,\infty[} + \frac{1}{2} 1_{[1+\alpha,\infty[}$ for $\alpha \in [0,1]$. The experiment is shown in figure \ref{fig:sep_2diracs}. There are 2 spectral separation for $\alpha = 0$. For $\alpha = 0.7$ and $\alpha = 1$ there is no more spectral separation.
	\item An exponentially-weighted distribution. The experiment is shown in figure \ref{fig:sep_exp}. There are 2 spectral separation for $\alpha = 0$, only 1 for $\alpha=2$ and no more for $\alpha = 5$.
\end{itemize}

\subsection{New spectral gaps due to weight gaps}
Silverstein and Choi \cite{Silverstein1995c} proved, in the case of the non-weighted sample covariance, that when $H$ is a finite sum of $K$ diracs, then there are at most $K-1$ spectral gaps in the asymptotic spectrum $F$.

In the case of the weighted covariance spectrum, this result is not true anymore, as gaps in $D$ can lead to new spectral gaps in $F$, sa explained theoretically in \cite{Couillet2015}.

To illustrate this new phenomenon, we consider:
\begin{itemize}
	\item $c = 0.05$ and $c=0.01$, 
	\item $H = 0.2 \times 1_{[1,\infty[} + 0.4 \times 1_{[3,\infty[} + 0.4 \times 1_{[10,\infty[}$,
	\item $D = \left(1 - \frac{1}{2\alpha} \right) \times 1_{[\frac{\alpha}{2\alpha - 1},\infty[} + \frac{1}{2\alpha} 1_{[\alpha,\infty[}$ for $\alpha \in \{1, 50 \}$ with $c = 0.05$ and $\alpha \in \{1, 200 \}$ with $c = 0.01$.
\end{itemize}

In  the standard non-weighted sample covariance case, we expect to find two spectral gaps roughly when $c<0.34$, one spectral gap for $c \in [0.34,0.40]$ and no spectral gaps for $c>0.40$. These values result from the theoretical study of Silverstein and Choi \cite{Silverstein1995c}.

\begin{figure}[ht]
\centering
\includegraphics[width=0.7\linewidth]{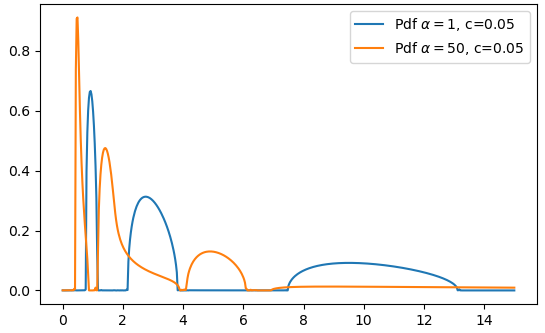}
\includegraphics[width=0.7\linewidth]{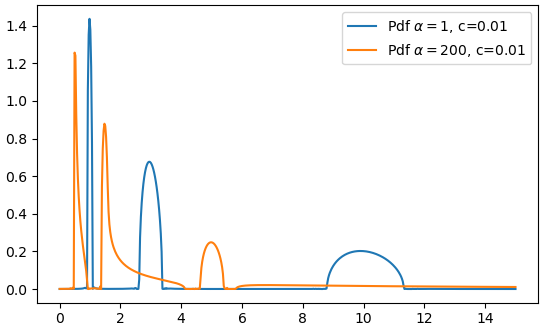}
\caption{Approximations of the asymptotic density of weighted sample eigenvalues for $c=0.05$ and $c=0.01$, with $H = 0.2 \times 1_{[1,\infty[} + 0.4 \times 1_{[3,\infty[} + 0.4 \times 1_{[10,\infty[}$ and weight distribution from a mixture of 2 diracs of parameter $\alpha$ for spectral gaps.}
\label{fig:sep_3}
\end{figure}

In our setting, we indeed have at most two spectral gaps with $\alpha = 1$, which coincides with the non-weighted situation. However, for $\alpha = 50$ and $\alpha = 200$, a third spectral gap appears at smaller values of $c$, respectively numerically for $c\leq 0.05$ and $c\leq 0.01$. Results are shown in figure \ref{fig:sep_3}. 

\subsection{Speed of convergence with heavy tails}
This last experimental part aims at empirically showing how the difference in finite sample to the asymptotic density evolves depending on the heaviness of tails.

To visualize this point, we are plotting the histogram of the empirical spectrum of the weighted sample covariance along with the theoretical asymptotic density. 

We draw the $(X_{ij})$ from a normalized t-distribution with varying degree of freedom $\nu \in \{2.5, 3.5, 4\}$. 

\begin{figure}[]
\centering
\includegraphics[width=0.7\linewidth]{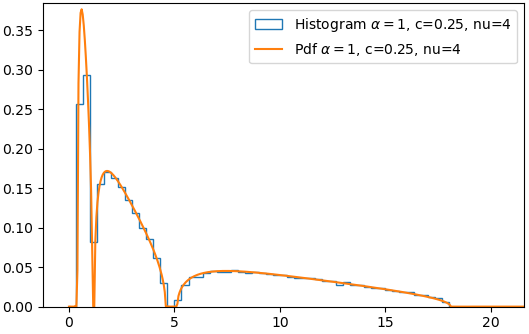}
\includegraphics[width=0.7\linewidth]{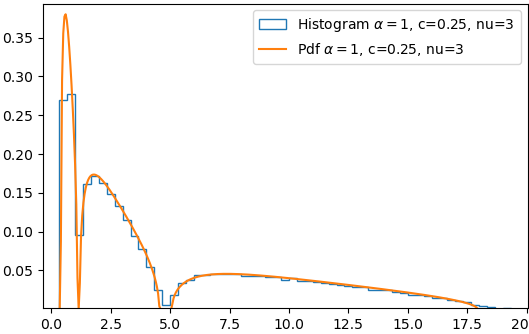}
\includegraphics[width=0.7\linewidth]{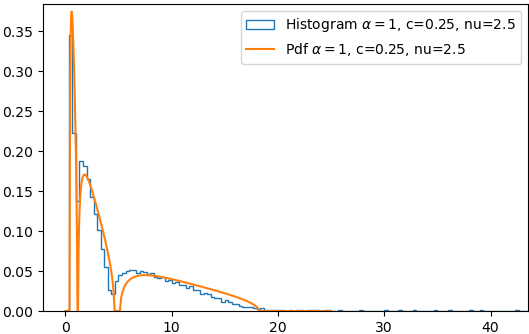}
\caption{Histogram from a t-distribution with $\nu = 4, \nu = 3$ and $\nu=2.5$ from top to bottom, and asymptotic theoretical density, for $n = 3000$, $c=0.25$, $H = 0.2 \times 1_{[1,\infty[} + 0.4 \times 1_{[3,\infty[} + 0.4 \times 1_{[10,\infty[}$ and exponentially-weighted distribution of parameter $\alpha = 1$.}
\label{fig:hist_exp}
\end{figure}

As previously, we chose $H = 0.2 \times 1_{[1,\infty[} + 0.4 \times 1_{[3,\infty[} + 0.4 \times 1_{[10,\infty[}$ and an exponentially-weighted distribution of parameter $\alpha = 1$ for the weights. The concentration ratio is set at $c=0.25$ and the dimension at $n = 3000$. The convergence is meant to be almost sure, so we did only one draw of spectrum per plot. The results are shown in figure \ref{fig:hist_exp}.

The conclusion of this experiment is that heavy tails tend to create rare but high eigenvalues, making the convergence slower as the tail grows. We recall that the convergence of the distributions is shown under weak convergence. So, the presence of very high eigenvalues far from the support of $F$ is not a contradiction as long as their frequence of appearance converges to zero, which experimentally is the case. This type of behavior is clearly visible for $\nu = 2.5$ in Figure \ref{fig:hist_exp} where eigenvalues higher than $40$ are sampled, but in very low proportion.

\section{Conclusion}
In this work, we provide a detailed introduction to the high dimensional spectrum of weighted sample covariance. We expose the setting, the convergence result in high dimension, along with a proof arguably simpler and more concise than the more optimal proof of \cite{Zhang2007}, and we illustrate some basic phenomenons to give the reader the correct intuition on what is happening in high dimensions when weights get involved in the estimation.

\section*{Acknowledgment}
We would like to thank Alexandre Miot and Gabriel Turinici for their insights and advices all along this work.

\newpage

\section{Appendix - Proof of Theorem \ref{dFE}}
\subsection{Truncation and centralization} 
In this first part of the proof, we use a similar approach as \cite{Silverstein1995b}, aiming at truncating and centering $Z$, $T$ and $W$ while conserving the same asymptotic spectrum for $F_n$.

\subsubsection{Truncation of $Z$ by $\sqrt{N}$}
Let $\hat Z_{ij} = Z_{ij} \mathrm{1}_{|Z_{ij}| < \sqrt{N}}$ and $\hat{\underline{B}}_n = \frac{1}{N} W^{1/2} \hat{Z}^* T \hat{Z} W^{1/2}$. We use Lemmas 2.5.a and 2.1.d \cite{Silverstein1995b}, and with $\lVert \cdot \rVert$ denoting the $\infty$-norm on bounded real functions, we have:
\begin{equation}\label{}
\begin{aligned}
	\lVert F^{\underline{B}_n} - F^{\hat{\underline{B}}_n} \rVert &\leq \frac{2}{N} \rk\left(W^{1/2}Z - W^{1/2} \hat Z \right) \\
	&\leq \frac{2}{N} \rk\left(Z - \hat Z \right) \\
	\lVert F^{\underline{B}_n} - F^{\hat{\underline{B}}_n} \rVert &\leq \frac{2}{N} \underbrace{\left|\left \{ i \in \llbracket 1,n \rrbracket, j \in \llbracket 1,N \rrbracket | Z_{ij} - \hat{Z}_{ij} \neq 0 \right \} \right|}_{\xi_n}
\end{aligned}
\end{equation} 
At this point, we follow the idea of the proof p.58-59 of \cite{Yin1986}. We have that:
\begin{equation}\label{}
\begin{aligned}
	\xi_n = \sum_{i,j} Y_{i,j} \text{ with } Y_{ij} = \mathrm{1}_{|Z_{ij}| \geq \sqrt{N}}.
\end{aligned}
\end{equation} 
Remark that the $Y_{ij}$ are i.i.d Bernoulli $\mathcal{B}(\eta)$ with $\eta = \p\left(|Z_{11}| \geq \sqrt{N} \right)$. As $\E[|Z_{11}|^2] = 1$, $\eta = o\left(\frac{1}{n}\right)$. 

Let $\epsilon > 0$. For $n$ large enough, we have:
\begin{equation}\label{}
\begin{aligned}
	\p\left(\frac{2}{N} \xi_n \geq \epsilon \right) &= \p\left(\frac{1}{nN} \sum_{ij} Y_{ij} \geq \frac{\epsilon}{2n} \right)  \leq \p\left(\frac{1}{nN} \sum_{ij} Y_{ij}\geq \eta + \frac{\epsilon}{4n} \right).
\end{aligned}
\end{equation}
As $\E[Y_{ij}] = \eta$, we have by Chernoff-Hoeffding's theorem, for $n$ large enough:
\begin{equation}\label{}
\begin{aligned}
	\p\left(\frac{2}{N} \xi_n \geq \epsilon \right) & \leq \exp\left[nN\left(\left(\eta + \frac{\epsilon}{4n} \right) \ln{\frac{\eta}{\eta + \frac{\epsilon}{4n}}} + \left(1 - \eta - \frac{\epsilon}{4n} \right) \ln{\frac{1-\eta}{1-\eta - \frac{\epsilon}{4n}}} \right) \right] \\
	& \leq \exp\left[nN\left(- \frac{\left(\frac{\epsilon}{4n} \right)^2}{2\left(\eta + \frac{\epsilon}{4n} \right)} - \frac{\left(\frac{\epsilon}{4n} \right)^3}{3\left(\eta + \frac{\epsilon}{4n} \right)^2} - \frac{ \left(\frac{\epsilon}{4n} \right)^2}{ 2\left(1 - \eta - \frac{\epsilon}{4n} \right)} + \frac{ \left(\frac{\epsilon}{4n} \right)^3}{3 \left(1 - \eta - \frac{\epsilon}{4n} \right)^2} \right) \right] \\
	&= \exp\left[- \frac{\epsilon^2N}{32n}\times \frac{1}{\left(\eta + \frac{\epsilon}{4n} \right)\left(1 - \eta - \frac{\epsilon}{4n} \right)} - \frac{\epsilon^3N}{3 \times 4^3n^2}\times \frac{1 - 2\left(\eta + \frac{\epsilon}{4n}\right)}{\left(\eta + \frac{\epsilon}{4n} \right)^2\left(1 - \eta - \frac{\epsilon}{4n} \right)^2}\right] \\
	\p\left(\frac{2}{N} \xi_n \geq \epsilon \right)& \leq \exp\left[- \epsilon N \times \left(\frac{1}{16} + \frac{2}{3 \times 4^3} \right) \right].
\end{aligned}
\end{equation}
Consequently, 
\begin{equation}\label{}
\begin{aligned}
	\lVert F^{\underline{B}_n} - F^{\hat{\underline{B}}_n} \rVert \underset{n \rightarrow \infty}{\longrightarrow} 0 \text{ almost surely.}
\end{aligned}
\end{equation}
\subsubsection{Centralization of $\hat Z$}
Let $\tilde Z_{ij} = \hat Z_{ij} - \E\left[\hat Z_{ij} \right]$ and $\tilde{\underline{B}}_n = \frac{1}{N} W^{1/2} \tilde{Z}^* T \hat{Z} W^{1/2}$. The $\hat Z_{ij}$ are i.i.d, so $\forall i,j, \E\left[\hat Z_{ij} \right] = \E\left[\hat Z_{11} \right]$. Consequently,
\begin{equation}\label{}
\begin{aligned}
	&\rk\left(W^{1/2} \tilde Z - W^{1/2} \hat Z \right) = \rk\left(W^{1/2} \E[\hat Z] \right) \leq 1.
\end{aligned}
\end{equation}
Using Lemma 2.5.a \cite{Silverstein1995b}, we have:
\begin{equation}\label{}
\begin{aligned}
	\lVert F^{\tilde{\underline{B}}_n} - F^{\hat{\underline{B}}_n} \rVert \leq \frac{2}{N} \rk\left(W^{1/2} \tilde Z - W^{1/2} \hat Z \right)  \leq \frac{2}{N} \underset{n \rightarrow \infty}{\longrightarrow} 0.
\end{aligned}
\end{equation}

\subsubsection{Truncation of $T$ and $W$ by $\alpha_n$ and $\beta_n$ respectively}
We denote $U$ a matrix of eigenvector of $T$ with associated eigenvalues $(\tau_1,...,\tau_n)$, so that $T = U \Diag\left((\tau_i)\right) U^*$. For $\alpha > 0$, we define:
\begin{equation}\label{}
\begin{aligned}
	T_\alpha = U \Diag\left(\left(\tau_i \mathrm{1}_{|\tau_i|\leq \alpha}\right)\right) U^*.
\end{aligned}
\end{equation}
If $\alpha$ and $-\alpha$ are continuity points of $H$, using Lemma 2.5.b \cite{Silverstein1995b}, we have that for any $n \times N$ Hermitian matrix $Q$:
\begin{equation}\label{}
\begin{aligned}
	\lVert F^{Q^* T Q} - F^{Q^* T_\alpha Q} \rVert \leq \frac{1}{N} \rk(T - T_\alpha) = \frac{1}{N} \sum_{i=1}^n \mathrm{1}_{|\tau_1|>\alpha}.
\end{aligned}
\end{equation}
Using Assumptions (b) and (c), we have then:
\begin{equation}\label{}
\begin{aligned}
	\frac{1}{N} \rk(T - T_\alpha) \underset{n \rightarrow \infty}{\longrightarrow} cH([-\alpha,\alpha]^c) \text{ almost surely.}
\end{aligned}
\end{equation}
Consequently,
\begin{equation}\label{}
\begin{aligned}
	\alpha := \alpha_n \underset{n \rightarrow \infty}{\longrightarrow} +\infty \implies \lVert F^{Q^* T Q} - F^{Q^* T_{\alpha} Q} \rVert \underset{n \rightarrow \infty}{\longrightarrow} 0 \text{ almost surely.}
\end{aligned}
\end{equation}
Similarly, for $\beta > 0$, we define:
\begin{equation}\label{}
\begin{aligned}
	W_\beta = \Diag\left(\left(W_{ii} \mathrm{1}_{|W_{ii}|\leq \beta}\right)\right).
\end{aligned}
\end{equation}
Using Lemma 2.5.b \cite{Silverstein1995b} and Assumptions (b) and (c), we have also for any $n \times N$ Hermitian matrix $Q$:
\begin{equation}\label{}
\begin{aligned}
	\lVert F^{W^{1/2} Q^* T_\alpha Q W^{1/2}} - F^{W_\beta^{1/2} Q^* T_\alpha Q W_\beta^{1/2}} \rVert \leq \frac{1}{n} \rk(W - W_\beta) \underset{n \rightarrow \infty}{\longrightarrow} \frac{1}{c}D([-\beta,\beta]^c) \text{ almost surely.}
\end{aligned}
\end{equation}
Consequently,
\begin{equation}\label{}
\begin{aligned}
	\beta := \beta_n \underset{n \rightarrow \infty}{\longrightarrow} +\infty \implies \lVert F^{W^{1/2} Q^* T_\alpha Q W^{1/2}} - F^{W_\beta^{1/2} Q^* T_\alpha Q W_\beta^{1/2}} \rVert \underset{n \rightarrow \infty}{\longrightarrow} 0 \text{ almost surely.}
\end{aligned}
\end{equation}
In the following, we choose $\alpha_n \uparrow +\infty$ and $\beta_n \uparrow +\infty$ so that:
\begin{itemize}
		\item $ (\alpha_n\beta_n)^4\left( \E\left[|Z_{11}|^2 \mathrm{1}_{|Z_{11}| \geq \ln N} \right] + \frac{1}{N}\right) \underset{n \rightarrow \infty}{\longrightarrow} 0, $
		\item $ \sum_{n=1}^\infty \frac{(\beta_n\alpha_n)^8}{N^2} \left( \E\left[|Z_{11}|^4 \mathrm{1}_{|Z_{11}| \geq \sqrt{N}} \right] +1  \right) < +\infty. $
\end{itemize}
We denote $\tilde{\underline{B}}_{\alpha,\beta} = \frac{1}{N}W_\beta^{1/2} \tilde{Z}^* T_\alpha \tilde{Z} W_\beta^{1/2}$. With $\alpha, \beta$ chosen following the rules above, we have:
\begin{equation}\label{}
\begin{aligned}
	\lVert F^{\tilde{\underline{B}}_{\alpha,\beta}} - F^{\tilde{\underline{B}}_{n}} \rVert \underset{n \rightarrow \infty}{\longrightarrow} 0 \text{ almost surely.}
\end{aligned}
\end{equation}

\subsubsection{Truncation of $\tilde Z$ by $\ln N$}
Let $\bar Z_{ij} = \tilde Z_{ij} \mathrm{1}_{|Z_{ij}| < \ln N} - \E\left[\tilde Z_{ij} \mathrm{1}_{|Z_{ij}| < \ln N} \right]$ and $\bar{\bar Z}_{ij} = \tilde Z_{ij} - \bar Z_{ij}$. Remark that $\E[\bar{\bar Z}] = \E[\bar Z] = 0$. Moreover, since $|\tilde Z_{ij}| \leq \ln N + \E[|Z_{11}|]$, we have for $n$ sufficiently large and some $a > 2$:
\begin{equation}\label{}
\begin{aligned}
	\frac{|\bar Z_{ij}|}{\sqrt{\E[|\bar Z_{11}|^2]}} \leq \frac{\ln N + \E[|Z_{11}|]}{\sqrt{\E[|\bar Z_{11}|^2]}} \leq a \ln n := \log n.
\end{aligned}
\end{equation}
We denote $\log n$ the logarithm of $n$ in base $\exp(1/a)$.

In this part, we use a metric $D$ on $\mathcal{M}(\R)$, the set of all subprobability distribution functions on $\R$, defined in \cite{Silverstein1995b} p.191. 
\begin{dft}[$D$-metric, \cite{Silverstein1995b}]
	Let $\{f_i \}$ be an enumeration of all continuous functions that take a constant $1/m$ value ($m$ a positive integer) on $[a,b]$, where $a, b$ are rational, $0$ on $]-\infty, a-1/m] \cup [b+1/m, +\infty[$, and linear on each $[a-1/m,a]$, $[b,b+1/m]$. For $F_1, F_2 \in \mathcal{M}(\R)$, we define:
	\begin{equation}\label{}
	\begin{aligned}
		D(F_1,F_2) := \sum_{i=1}^\infty \left|\int f_i dF_1 - \int f_idF_2 \right|2^{-i}.
	\end{aligned}
	\end{equation}
	$D$ is a metric on $\mathcal{M}(\R)$ inducing the topology of vague convergence.
\end{dft}
We denote $\bar{\underline{B}}_{\alpha,\beta} = \frac{1}{N} W_\beta^{1/2} \bar{Z}^* T_\alpha \bar{Z} W_\beta^{1/2}$. Using Lemma 2.1.c and Equation (2.4) \cite{Silverstein1995b}, we have:
\begin{equation}\label{}
\begin{aligned}
	D^2\left(F^{\tilde{\underline{B}}_{\alpha,\beta}}, F^{\bar{\underline{B}}_{\alpha,\beta}} \right) 
	& \leq &&\frac{1}{N} \tr \left[\left( W_\beta^{1/2} \tilde{Z}^* T_\alpha \tilde{Z} W_\beta^{1/2} - W_\beta^{1/2} \bar{Z}^* T_\alpha \bar{Z} W_\beta^{1/2}\right)^2 \right] \\
	& = &&\frac{1}{N} \tr \left[\left( W_\beta^{1/2} (\bar{\bar Z} + \bar Z)^* T_\alpha (\bar{\bar Z} + \bar Z) W_\beta^{1/2} - W_\beta^{1/2} \bar{Z}^* T_\alpha \bar{Z} W_\beta^{1/2}\right)^2 \right] \\
	D^2\left(F^{\tilde{\underline{B}}_{\alpha,\beta}}, F^{\bar{\underline{B}}_{\alpha,\beta}} \right) 
	& = &&\frac{1}{N} \Bigg( \tr \left[\left( W_\beta^{1/2} \bar{\bar Z} ^* T_\alpha \bar{\bar Z} W_\beta^{1/2}\right)^2\right] 
			+ \tr\left[\left(W_\beta^{1/2} \bar{\bar Z} ^* T_\alpha \bar Z W_\beta^{1/2} + W_\beta^{1/2} \bar Z ^* T_\alpha \bar{\bar Z} W_\beta^{1/2} \right)^2 \right] \\
	&    &&	+ 2\tr\left[\left(W_\beta^{1/2} \bar{\bar Z} ^* T_\alpha \bar{\bar Z} W_\beta^{1/2} \right)\left( W_\beta^{1/2} \bar{\bar Z} ^* T_\alpha \bar Z W_\beta^{1/2} + W_\beta^{1/2} \bar Z ^* T_\alpha \bar{\bar Z} W_\beta^{1/2}\right) \right]  \Bigg).
\end{aligned}
\end{equation}
Using Cauchy-Schwarz inequality, we have:
\begin{equation}\label{}
\begin{aligned}
	D^2\left(F^{\tilde{\underline{B}}_{\alpha,\beta}}, F^{\bar{\underline{B}}_{\alpha,\beta}} \right) 
	& = &&\frac{1}{N} \Bigg( \tr \left[\left( W_\beta^{1/2} \bar{\bar Z} ^* T_\alpha \bar{\bar Z} W_\beta^{1/2}\right)^2\right] 
			+ 4\tr\left[W_\beta^{1/2} \bar{\bar Z} ^* T_\alpha \bar Z W_\beta \bar Z ^* T_\alpha \bar{\bar Z} W_\beta^{1/2} \right] \\
	&    &&	+ 4\sqrt{\tr \left[\left( W_\beta^{1/2} \bar{\bar Z} ^* T_\alpha \bar{\bar Z} W_\beta^{1/2}\right)^2\right] \tr\left[W_\beta^{1/2} \bar{\bar Z} ^* T_\alpha \bar Z W_\beta \bar Z ^* T_\alpha \bar{\bar Z} W_\beta^{1/2}  \right]}  \Bigg).
\end{aligned}
\end{equation}
We show, through Von Neumann's trace inequality, that:
\begin{equation}\label{}
\begin{aligned}
	\tr \left[\left( W_\beta^{1/2} \bar{\bar Z} ^* T_\alpha \bar{\bar Z} W_\beta^{1/2}\right)^2\right] &\leq \alpha^2\beta^2 \tr\left[\left(\bar{\bar Z}^*\bar{\bar Z} \right)^2 \right], \\
	\text{and }\tr\left[W_\beta^{1/2} \bar{\bar Z} ^* T_\alpha \bar Z W_\beta \bar Z ^* T_\alpha \bar{\bar Z} W_\beta^{1/2} \right]  &\leq \alpha^2 \beta ^2 \sqrt{\tr\left[\left(\bar{\bar Z}^*\bar{\bar Z} \right)^2 \right]\tr\left[\left({\bar Z}^*{\bar Z} \right)^2 \right]}.
\end{aligned}
\end{equation}
And, from \cite{Silverstein1995b} p.185 and Assumption (b), we have:
\begin{equation}\label{}
\begin{aligned}
	\frac{1}{N}\tr\left[\left({\bar Z}^*{\bar Z} \right)^2 \right]  \underset{n \rightarrow \infty}{\longrightarrow} c(1+c) \text{ almost surely}.
\end{aligned}
\end{equation}
We have also from \cite{Silverstein1995b} p.185, with $K$ and $K'$ constant independent of $n$,
\begin{equation}\label{}
\begin{aligned}
	&\E\left[\frac{1}{N}\tr\left[\left(\bar{\bar Z}^*\bar{\bar Z}\right)^2 \right]\right]  \leq K'\left(\E\left[|Z_{11}|^2 \mathrm{1}_{|Z_{11}| \geq \ln N} \right] + \frac{1}{N} \right), \\
	&\V\left[\frac{1}{N}\tr\left[\left(\bar{\bar Z}^*\bar{\bar Z}\right)^2 \right]\right]  \leq \frac{K}{N^2} \left(  \E\left[|Z_{11}|^4 \mathrm{1}_{|Z_{11}| \geq \sqrt{N}} \right] +1 \right).
\end{aligned}
\end{equation}
So, $\E\left[\frac{\alpha^4 \beta^4}{N}\tr\left[\left(\bar{\bar Z}^*\bar{\bar Z}\right)^2 \right]\right] \underset{n \rightarrow \infty}{\longrightarrow} 0$ and $\V\left[\frac{\alpha^4 \beta^4}{N}\tr\left[\left(\bar{\bar Z}^*\bar{\bar Z}\right)^2 \right]\right]$ is summable, so:
\begin{equation}\label{}
\begin{aligned}
	\frac{\alpha^4 \beta^4}{N}\tr\left[\left(\bar{\bar Z}^*\bar{\bar Z}\right)^2 \right] \underset{n \rightarrow \infty}{\longrightarrow} 0 \text{ almost surely}.
\end{aligned}
\end{equation}
Backing up, we have then:
\begin{equation}\label{}
\begin{aligned}
	D^2\left(F^{\tilde{\underline{B}}_{\alpha,\beta}}, F^{\bar{\underline{B}}_{\alpha,\beta}} \right)  \underset{n \rightarrow \infty}{\longrightarrow} 0 \text{ almost surely}.
\end{aligned}
\end{equation}
So, $F^{\tilde{\underline{B}}_{\alpha,\beta}} - F^{\bar{\underline{B}}_{\alpha,\beta}} \overset{v}{\longrightarrow} 0 $, $ \overset{v}{\longrightarrow}$ denotes the vague convergence. 

\subsubsection{Conclusion of the truncations and centralizations}
The conclusion of that first part of the proof is that it is sufficient to show that $F^{\bar{\underline{B}}_{\alpha,\beta}}  \overset{v}{\longrightarrow} \underline{F}$ for some $ \underline{F} \in \mathcal{M}(\R)$ (we will in fact prove the weak convergence) in order to prove that $F^{ \underline{B}_n} \overset{v}{\longrightarrow}  \underline{F}$, which is equivalent to show that $F^{{B}_n} \overset{v}{\longrightarrow} F := \frac{1}{c} \underline{F} + \frac{1-c}{c} \mathrm{1}_{[0,\infty[}$.

Moreover, as by definition $F^{B_n}(\R) = 1$, if $F(\R) = 1$, it is equivalent to prove that $F^{B_n} \underset{n \rightarrow \infty}{\implies} F$. 

In the following, we focus on the truncated and centralized variables. In fact, we swap:
\begin{itemize}
	\item $Z_{ij}$ by $\bar Z_{ij}/\sqrt{\E[|\bar Z_{ij}|^2]}$,
	\item $T$ by $\E[|\bar Z_{11}|^2] T_\alpha$,
	\item $W$ by $W_\beta$.
\end{itemize}
Remark that we can impose $\alpha \leq \log n$ and $\beta \leq \log n$ without compromising any of the properties on $\alpha$ and $\beta$ defined in Section 6.1.3.

In addition to (a)-(e), we can now use the following Assumptions.
\begin{as}
	Assume that:
	\begin{itemize}
		\item[(f)] $\forall i,j, |Z_{ij}| \leq \log n$,
		\item[(g)] $\lVert T \rVert \leq \log n$,
		\item[(h)] $\lVert W \rVert \leq \log n$.
	\end{itemize}
\end{as}

From now, $B_n = \frac{1}{N} T^{1/2} Z^* W Z T^{1/2}$ and $\underline{B}_n = \frac{1}{N} W^{1/2} Z T Z^* W^{1/2}$ use the truncated variables $Z$, $T$ and $W$, and so are $F_n$, $\underline{F}_n$, $m_n$ and $\underline{m}_n$.

The following aims at proving, with those truncated variables, that for some p.d.f $F$, $F^{B_n} \overset{v}{\longrightarrow} F$ a.s., which proves the theorem. For that, we prove that for some p.d.f $F$, $F^{B_n} \underset{n \rightarrow \infty}{\implies} F$ a.s. It is done using the Cauchy-Stieltjes transform: we prove that for all $z \in \C_+$, a.s. $m_n(z) \underset{n \rightarrow \infty}{\longrightarrow} m(z)$ where $m$ is the Cauchy-Stieltjes transform of a p.d.f.

We introduce new objects of interest in the analysis, namely $q_j, r_j$, and $B_{(j)}$.
\begin{notat}
	For $j \in \llbracket 1,N \rrbracket$, we denote:
	\begin{itemize}
		\item $q_j = \frac{1}{\sqrt{n}}Z_{\cdot j}$, 
		\item $r_j =  \frac{1}{\sqrt{N}} T^{1/2}Z_{\cdot j} W_{jj}^{1/2}$, 
		\item $B_{(j)} = B_n - r_j r_j^*$, 
	\end{itemize}
\end{notat}

\subsection{Concentration}
\subsubsection{Preliminary derivations}
Let's introduce a technical result.
\begin{lem}[Eq (2.1) \cite{Silverstein1995c}]\label{tech0}
	For $B$ a $(n,n)$ matrix, $q \in  \C^n$ for which $B$ and $B + qq^*$ is invertible, we have:
	\begin{equation}\label{}
	\begin{aligned}
		q^* (B + qq^*)^{-1} = \frac{1}{1 + q^*B^{-1}q}q^* B^{-1}.
	\end{aligned}
	\end{equation}
\end{lem}

We have: 
\begin{equation}\label{}
\begin{aligned}
	(B_n -zI) +zI &= \sum_{j=1}^N r_j r_j^*  \\
	I + z(B_n - zI)^{-1} &= \sum_{j=1}^N r_j r_j^* (B_n -zI)^{-1} \\
	I + z(B_n - zI)^{-1} &= \sum_{j=1}^N \frac{1}{1 + r_j^* (B_{(j)} -zI)^{-1}r_j} r_j r_j^* (B_{(j)} -zI)^{-1}  \\
	c_n + z c_n m_n(z) &= 1 - \frac{1}{N} \sum_{j=1}^N \frac{1}{1 + r_j^* (B_{(j)} -zI)^{-1}r_j} \\
	\underline{m}_n(z) &= - \frac{1}{N} \sum_{j=1}^N \frac{1}{z\left(1 + r_j^* (B_{(j)} -zI)^{-1}r_j\right)}.
\end{aligned}
\end{equation} 

Here is the crucial difference between the proof in \cite{Silverstein1995a} in the evenly weighted case and the weighted case. We denote:
\begin{equation}\label{} 
\begin{aligned}
	\alpha_n = - \frac{1}{z\underline{m}_n(z)} \frac{1}{N} \sum_{j=1}^N  W_{jj}\left(r_j^* (B_n -zI)^{-1}r_j - 1 \right),
\end{aligned}
\end{equation} 
while in \cite{Silverstein1995c}, $\alpha_n = 1$ is used.

We have equivalently: 
\begin{equation}\label{}
\begin{aligned}
	\alpha_n &= - \frac{1}{z\underline{m}_n(z)} \frac{1}{N} \sum_{j=1}^N  \frac{W_{jj}}{1 + r_j^* (B_{(j)} -zI)^{-1}r_j} \\
	\text{and } \alpha_n &= \frac{\frac{1}{N} \sum_{j=1}^N \frac{W_{jj}}{1 + r_j^* (B_{(j)} -zI)^{-1}r_j}}{\frac{1}{N} \sum_{j=1}^N \frac{1}{1 + r_j^* (B_{(j)} -zI)^{-1}r_j}},
\end{aligned}
\end{equation} 

Then:
\begin{equation}\label{}
\begin{aligned}
	&(-z \alpha_n \underline{m}_n(z) T_n - zI)^{-1} - (B_n - zI)^{-1} = \\
	&(-z \alpha_n \underline{m}_n(z) T_n - zI)^{-1}\left(z \alpha_n \underline{m}_n(z) T_n + \sum_{j=1}^N r_j r_j^*  \right) (B_n - zI)^{-1} =  \\
	&(-z \alpha_n \underline{m}_n(z) T_n - zI)^{-1}\sum_{j=1}^N \frac{r_j r_j^* (B_{(j)} -zI)^{-1} - \frac{1}{N}\alpha_n T_n (B_n -zI)^{-1}}{1 + r_j^* (B_{(j)} -zI)^{-1}r_j}.
\end{aligned}
\end{equation} 
Applying the trace and dividing by $n$, we obtain:
\begin{equation}\label{}
\begin{aligned}
	&\frac{1}{n} \tr\left((-z \alpha_n \underline{m}_n(z) T_n - zI)^{-1}\right) - m_n(z) = &&\\
	&-\frac{1}{N}\sum_{j=1}^N \frac{1}{z(1 + r_j^* (B_{(j)} -zI)^{-1}r_j)} \Big(W_{jj}q_j^* T^{1/2} (B_{(j)} - zI)^{-1} (\alpha_n \underline{m}_n(z) T_n + I)^{-1}T^{1/2} q_j \\
	& \qquad - \frac{1}{n} \tr \left( (\alpha_n \underline{m}_n(z) T_n + I)^{-1} \alpha_n T_n (B_n - zI)^{-1}\right)\Big)\\
\end{aligned}
\end{equation} 
So,
\begin{equation}\label{eq_cvg}
\begin{aligned}
	&\frac{1}{n} \tr\left((-z \alpha_n \underline{m}_n(z) T_n - zI)^{-1}\right) - m_n(z) = &&\\
	&-\frac{1}{zN}\sum_{j=1}^N q_j^* \left(\frac{W_{jj}T_n^{1/2} (B_{(j)} - zI)^{-1} (\alpha_n \underline{m}_n(z) T_n + I)^{-1}T_n^{1/2}}{1 + r_j^* (B_{(j)} -zI)^{-1}r_j} \right)q_j \\
	& + \frac{1}{zNn} \tr \left( (\alpha_n \underline{m}_n(z) T_n + I)^{-1} \sum_{j=1}^N \frac{\alpha_n T_n}{1 + r_j^* (B_{(j)} -zI)^{-1}r_j} (B_n - zI)^{-1}\right) = \\
	&-\frac{1}{N} \sum_{j=1}^N \frac{1}{z(1 + r_j^* (B_{(j)} -zI)^{-1}r_j)} d_j,
\end{aligned}
\end{equation} 
where:
\begin{equation}\label{}
\begin{aligned}
	&d_j = &&W_{jj}q_j^* T_n^{1/2} (B_{(j)} - zI)^{-1} (\alpha_n \underline{m}_n(z) T_n + I)^{-1}T_n^{1/2}q_j \\
	& &&- \frac{1}{n}\tr \left(W_{jj} (\alpha_n \underline{m}_n(z) T_n + I)^{-1}  T_n (B_n - zI)^{-1}\right).
\end{aligned}
\end{equation}

The strategy of the proof is the following:
\begin{itemize}
	\item prove that $\max_{j \leq N} d_j \longrightarrow 0$ a.s., and that $\frac{1}{n} \tr\left((-z \alpha_{n} \underline{m}_n(z) T_n - zI)^{-1}\right) - m_n(z) \longrightarrow 0$ a.s.,
	\item show that a.s. it exists $\tilde m(z) \in \C \backslash \C_+$ and a subsequence ${n_i}$ so that $\tilde m_{n_i}(z) := -z\alpha_{n_i}\underline{m}_{n_i}(z) \longrightarrow \tilde m(z)$ a.s.,
	\item prove that $\tilde m(z)$ is the unique solution in $\C_-$ of a functional equation, and deduce $\tilde m_{n}(z) \longrightarrow \tilde m(z)$ a.s.,
	\item deduce that a.s. it exists $ m(z) \in \C_+$, uniquely defined in function of $\tilde m(z)$ so that $m_n(z) \longrightarrow m(z)$ a.s.,
	\item similarly deduce that a.s. it exists $ \Theta^{(1)}(z) \in \C_+$, uniquely defined in function of $\tilde m(z)$ so that $\Theta^{(1)}_n(z) \longrightarrow \Theta^{(1)}(z)$ a.s.,
	\item conclude proving that  $m$ is the Cauchy-Stieltjes transform of a p.d.f.
\end{itemize}

\subsubsection{Decomposition of $d_j$}
Much of the truth of this proof relies upon the following lemma from \cite{Silverstein1995b}.
\begin{lem}[Lemma 3.1 \cite{Silverstein1995b}]\label{lemma_rmt}
	Let $C$ a Hermitian $n \times n$ matrix so that $\lVert C \rVert \leq 1$, and $Y = (Z_1,...,Z_n)^T$, $Z_i \in \C$ where the $Z_i$'s are independent, $\forall i, \E[Z_i] = 0$, $\E[|Z_i|^2] = 1$ and $|Z_i| \leq \log n$. Then,
	\begin{equation}\label{}
	\begin{aligned}
		&\E\left[\left| Y^* C Y - \tr C \right|^6 \right] \leq K n^3 \log(n)^{12},
	\end{aligned}
	\end{equation}
	where the constant $K$ does not depend on $n$, $C$, nor the distribution of $Z_1$.
\end{lem}

\begin{rmk}
	Lemma 3.1 \cite{Silverstein1995b} assumes that the $Z_i$ are i.i.d but the "indentically distributed" assumption is in fact never used in the proof.
\end{rmk}

In order to use this Lemma on $d_j$, we decompose it into negligible terms and a term of the form $q_j^* C q_j - \frac{1}{n}tr(C)$ where $C$ is independent of $q_j$. 

For that, we denote: 
\begin{equation}\label{}
\begin{aligned}
	\tilde m_{n}(z) &:= -z\alpha_n\underline{m}_{n}(z) \\
	\text{and } \tilde m_{(j)}(z) &:= \frac{1}{N} \sum_{i \neq j}  W_{ii}\left(r_i^* (B_{(j)} -zI)^{-1}r_i - 1 \right).
\end{aligned}
\end{equation}

We have the following decomposition of $d_j$: 
\begin{equation}\label{}
\begin{aligned}
	d_j &= d_j^{(1)} + d_j^{(2)} + d_j^{(3)} + d_j^{(4)},
\end{aligned}
\end{equation}
with:
\begin{equation}\label{}
\begin{aligned}
	d_j^{(1)} &= W_{jj}q_j^* T_n^{1/2} (B_{(j)} - zI)^{-1} \left[\left(-\frac{\tilde{m}_n(z)}{z}  T_n + I\right)^{-1} - \left(-\frac{\tilde{m}_{(j)}(z)}{z}  T_n + I\right)^{-1} \right]T_n^{1/2}q_j, \\
	d_j^{(2)} &= W_{jj}q_j^* T_n^{1/2} (B_{(j)} - zI)^{-1}\left(-\frac{\tilde{m}_{(j)}(z)}{z}  T_n + I\right)^{-1}T_n^{1/2}q_j \\
	& \qquad\qquad\qquad\qquad\qquad\qquad - \frac{W_{jj}}{n}\tr \left( \left(-\frac{ \tilde{m}_{(j)}(z)}{z} T_n + I\right)^{-1} T_n (B_{(j)} - zI)^{-1}\right), \\
	d_j^{(3)} &= \frac{W_{jj}}{n}\tr \left(  \left[\left(-\frac{\tilde{m}_{(j)}(z)}{z}  T_n + I\right)^{-1} - \left(-\frac{\tilde{m}_{n}(z)}{z}  T_n + I\right)^{-1} \right] T_n (B_{(j)} - zI)^{-1}\right), \\
	d_j^{(4)} &= \frac{W_{jj}}{n}\tr \left(\left(-\frac{\tilde{m}_{n}(z)}{z}  T_n + I\right)^{-1} T_n \left[(B_{(j)} - zI)^{-1} - (B_{n} - zI)^{-1}\right]\right).
\end{aligned}
\end{equation}

In order to prove that for each $k \in \llbracket 1, 4 \rrbracket, \max_{j \leq N} d_j^{(k)} \longrightarrow 0$ a.s., we need some technical lemmas. They essentially provide the necessary inequalities to prove that $d_j^{(1)}, d_j^{(3)}$ and $d_j^{(4)}$ are indeed negligible, to finally use Lemma \ref{lemma_rmt} on $d_j^{(2)}$ and prove that  $\max_{j \leq N} |d_j| \longrightarrow 0$ a.s.

\subsubsection{Technical lemmas}

\begin{lem}\label{ineq}
	We have the following inequalities, for $j \in \llbracket 1,N \rrbracket$:
	\begin{equation}\label{}
	\begin{aligned}
		&\lVert (B_n - zI)^{-1} \rVert \leq \frac{1}{v}, \\
		&\lVert (B_{(j)} - zI)^{-1} \rVert \leq \frac{1}{v}, \\
		&\frac{1}{|z(1 + r_j^* (B_{(j)} -zI)^{-1}r_j)|} \leq \frac{1}{v}.
	\end{aligned}
	\end{equation} 
\end{lem}

\begin{proof}
	Let $j \in \llbracket 1, N\rrbracket$. The two first inequalities comes from the fact that $B_n$ and $B_{(j)}$ are Hermitian, so for any eigenvalue $\lambda$ of $B_n - zI$ or $B_{(j)} - zI$, we have that $|\lambda| \geq |\Im[\lambda]| = v$. 
	
	For the third inequality, we remark that:
	\begin{equation}\label{}
	\begin{aligned}
		\Im r_j^* (B_{(j)}/z -I)^{-1} r_j & = \frac{1}{2i} r_j^* \left[(B_{(j)}/z -I)^{-1} - (B_{(j)}/z^* -I)^{-1} \right] r_j \\
		& = \frac{v}{|z|^2} r_j^* (B_{(j)}/z -I)^{-1}B_{(j)} (B_{(j)}/ z^* -I)^{-1} r_j  \\
		\Im r_j^* (B_{(j)}/z -I)^{-1} r_j & \geq 0.
	\end{aligned}
	\end{equation} 
	So, we deduce that:
	\begin{equation}\label{}
	\begin{aligned}
		\frac{1}{|z(1 + r_j^* (B_{(j)} -zI)^{-1}r_j)|} \leq \frac{1}{v}.
	\end{aligned}
	\end{equation} 
\end{proof}

\begin{lem}\label{tech1}
	We denote $\bar W =  \frac{1}{N} \sum_{i=1}^N W_{ii}$. For $z = u +iv$, $v >0$ and $j \in \llbracket 1, N\rrbracket$, we have for any nonnegative Hermitian matrix $A$:
	\begin{equation}\label{}
	\begin{aligned}
		&\left \lVert \left( - \frac{\tilde m_n(z)}{z} A + I \right)^{-1} \right \rVert \leq f(z,\lVert A \rVert), \text{ and }
		\left \lVert \left( - \frac{\tilde m_{(j)}(z)}{z} A + I \right)^{-1} \right \rVert \leq f(z,\lVert A \rVert),
	\end{aligned}
	\end{equation}
	where:
	\begin{equation}\label{}
	\begin{aligned}
		f(z,\lVert A \rVert) = \begin{cases}
			 \max \left(2, \frac{4}{v} \bar W\lVert A \rVert\right), &\text{ if } u= 0, \\
			 16 \left( \frac{|z|^2}{4v^2|u|}\bar W \lVert A \rVert +1 \right) \times \max \left(\frac{1}{3}, \frac{|u|}{v} \right), &\text{ otherwise.}
		\end{cases}
	\end{aligned}
	\end{equation}
\end{lem}

\begin{proof}
	Let z = u + iv$, v> 0, u \in \R$. For $j \in \llbracket 1,N \rrbracket$, we denote by $(u_{ij})_{i=1}^n$ a set of eigenvectors of $B_{(j)}$ with associated eigenvalues $(\lambda_{ij})_{i=1}^n$. We then derive the following formulation:
		\begin{equation}\label{}
		\begin{aligned}
			R := \Re\left[- \frac{\tilde m_n(z)}{z} \right] &= - \frac{1}{N} \sum_{j=1}^N \frac{W_{jj}}{\left|z + r_j^* \left(\frac{B_{(j)}}{z} - I \right)^{-1}r_j\right|^2} \left(u + \sum_{i=1}^n \frac{|r_j^*u_{ij}|^2}{|\lambda_{ij} - z|^2} (\lambda_{ij}u - |z|^2) \right),\\
			I :=\Im\left[- \frac{\tilde m_n(z)}{z} \right] &= \frac{1}{N} \sum_{j=1}^N \frac{W_{jj}}{\left|z + r_j^* \left(\frac{B_{(j)}}{z} - I \right)^{-1}r_j\right|^2} \left(v + \sum_{i=1}^n \frac{|r_j^*u_{ij}|^2}{|\lambda_{ij} - z|^2} \lambda_{ij}v \right) \geq 0.
		\end{aligned}
		\end{equation}
		Using Cauchy-Schwarz inequality, we deduce:
		\begin{equation}\label{}
		\begin{aligned}
			\frac{1}{N} \sum_{j=1}^N \frac{W_{jj}|z|^2\sum_{i=1}^n \frac{|r_j^*u_{ij}|^2}{|\lambda_{ij} - z|^2}}{\left|z + r_j^* \left(\frac{B_{(j)}}{z} - I \right)^{-1}r_j\right|^2}  \leq &\sqrt{\frac{1}{N} \sum_{j=1}^N \frac{W_{jj}\left( \sum_{i=1}^n \frac{|r_j^*u_{ij}|^2}{|\lambda_{ij} - z|^2}\right)^2}{\left|1 + r_j^* \left(B_{(j)} - zI \right)^{-1}r_j\right|^2} } \sqrt{\frac{1}{N} \sum_{j=1}^N \frac{W_{jj}}{\left|z + r_j^* \left(\frac{B_{(j)}}{z} - I \right)^{-1}r_j\right|^2}}  \\
			\frac{1}{N} \sum_{j=1}^N \frac{W_{jj}|z|^2\sum_{i=1}^n \frac{|r_j^*u_{ij}|^2}{|\lambda_{ij} - z|^2}}{\left|z + r_j^* \left(\frac{B_{(j)}}{z} - I \right)^{-1}r_j\right|^2}  \leq & \frac{|z|\sqrt{\bar W}}{v\sqrt{v}} \sqrt{I}.\\
		\end{aligned}
		\end{equation}
		So, combining both previous points, we have:
		\begin{equation}\label{sqrt_ineq}
		\begin{aligned}
			\left| R \right|& \leq \frac{|u|}{v} I + \frac{|z|\sqrt{\bar W}}{v\sqrt{v}} \sqrt{I}.
		\end{aligned}
		\end{equation}
		
		Suppose $u \neq 0$. We denote $K:=\frac{|z|\sqrt{\bar W}}{2v\sqrt{|u|}}$. We have then:
		\begin{equation}\label{}
		\begin{aligned}
			&\sqrt{\frac{v}{|u|}} \left(-  K + \sqrt{ K^2 + \left| R \right|}\right) \leq \sqrt{I}.
		\end{aligned}
		\end{equation}
		Now, let $x \geq 0$. Then:
		\begin{equation}\label{}
		\begin{aligned}
			\left|- \frac{\tilde m_n(z)}{z}x + 1 \right|^2 &= \left(Rx +1\right)^2 + I^2 x^2\\
			\left|- \frac{\tilde m_n(z)}{z}x + 1 \right|^2& \geq (-|R|x +1)^2  + \frac{v^2}{|u|^2} \left(-  K\sqrt{x} + \sqrt{ K^2x + \left| R \right|x}\right)^4.
		\end{aligned}
		\end{equation}
		We denote: $t := \sqrt{ K^2x +\left| R \right|x}\in \R_+$ . We have then:
		\begin{equation}\label{}
		\begin{aligned}
			\left|- \frac{\tilde m_n(z)}{z}x + 1 \right|^2 &\geq (t^2 - K^2x - 1)^2 + \frac{v^2}{|u|^2} (t - K\sqrt{x})^4.
		\end{aligned}
		\end{equation}
		We denote $a:=K\sqrt{x}$, $b:= \sqrt{K^2x +1}$ and we split the study of the right part of the previous equation between $[0,(a+b)/2]$ and $[(a+b)/2,+\infty[$. The lower bounds rely mainly on the fact that $b- a \geq \frac{1}{2b}$:
		\begin{itemize}
			\item Let $t \in [0,(a+b)/2]$. Then,
				\begin{equation}\label{}
				\begin{aligned}
					(t^2 - K^2x - 1)^2 + \frac{v^2}{|u|^2} (t - K\sqrt{x})^4 &\geq \left(\left(\frac{a+b}{2}\right)^2 - b^2\right)^2 \\
					& = \frac{(a + 3b)^2(b-a)^2}{16} \\
					& \geq \frac{9(b-a)^4}{16} \\
					(t^2 - K^2x - 1)^2 + \frac{v^2}{|u|^2} (t - K\sqrt{x})^4 &\geq \frac{9}{16^2b^4}.
				\end{aligned}
				\end{equation}
			\item Let $t \in [(a+b)/2, +\infty[$. Then,
				\begin{equation}\label{}
				\begin{aligned}
					(t^2 - K^2x - 1)^2 + \frac{v^2}{|u|^2} (t - K\sqrt{x})^4 &\geq  \frac{v^2}{16|u|^2}(b-a)^4 \\
					(t^2 - K^2x - 1)^2 + \frac{v^2}{|u|^2} (t - K\sqrt{x})^4 &\geq  \frac{v^2}{16^2|u|^2b^4}. \\
				\end{aligned}
				\end{equation}
		\end{itemize}
		Backing up, we have that:
		\begin{equation}\label{}
		\begin{aligned}
			&\left|- \frac{\tilde m_n(z)}{z}x + 1 \right|^2 \geq \frac{1}{16^2\left(K^2x +1 \right)^2}\times \min\left(\frac{v^2}{|u|^2}, 9 \right).
		\end{aligned}
		\end{equation}
		It finally leads to:
		\begin{equation}\label{}
		\begin{aligned}
			&\left \lVert \left( - \frac{\tilde m_n(z)}{z} A + I \right)^{-1} \right \rVert \leq 16 \left( \frac{|z|^2}{4v^2|u|}\bar W \lVert A \rVert +1 \right) \times \max \left(\frac{1}{3}, \frac{|u|}{v} \right),
		\end{aligned}
		\end{equation}
		which proves the first inequality when $u \neq 0$.
		
		Now suppose $u= 0$. From Equation \eqref{sqrt_ineq}, we have:
		\begin{equation}\label{sqrt_ineq}
		\begin{aligned}
			\left| R \right|& \leq \sqrt{\frac{\bar W}{v}} \sqrt{I}.
		\end{aligned}
		\end{equation}
		If $\bar W=0$, then:
		\begin{equation}\label{}
		\begin{aligned}
			&\left|- \frac{\tilde m_n(z)}{z}x + 1 \right|^2 =1 + x^2 I^2 \geq 1 \geq \frac{1}{4}.
		\end{aligned}
		\end{equation}
		Otherwise, with $y = Rx$ and $a = \frac{v^2}{\bar W^4x^2}$:
		\begin{equation}\label{}
		\begin{aligned}
			&\left|- \frac{\tilde m_n(z)}{z}x + 1 \right|^2 \geq (Rx + 1)^2 +x^2 \frac{R^4v^2}{\bar W^4} = ay^4 + (y+1)^2.
		\end{aligned}
		\end{equation}
		Splitting the study between $]-\infty, 1/2]$ and $[1/2, +\infty[$, we have that $ay^4 + (y+1)^2 \geq \min \left(\frac{a}{16}, \frac{1}{4} \right).$ So, we deduce the first inequality when $u = 0$:
		\begin{equation}\label{}
		\begin{aligned}
			&\left \lVert \left( - \frac{\tilde m_n(z)}{z} A + I \right)^{-1} \right \rVert \leq \max\left(2, \frac{4\bar W\lVert A \rVert}{v} \right).
		\end{aligned}
		\end{equation}
		
		Let $j \in \llbracket 1, N \rrbracket$. Using the same method for $\left \lVert \left( - \frac{\tilde m_{(j)}(z)}{z} A + I \right)^{-1} \right \rVert$, we have with $\bar W_{(j)} := \frac{1}{N} \sum_{i \neq j} W_{ii}$ :
		\begin{equation}\label{}
		\begin{aligned}
			\left| \Re\left[- \frac{\tilde m_{(j)}(z)}{z} \right] \right|& \leq \frac{|u|}{v} \Im\left[- \frac{\tilde m_{(j)}(z)}{z} \right] + \frac{|z|\sqrt{\bar W_{(j)}}}{v\sqrt{v}}  \sqrt{ \Im\left[- \frac{\tilde m_{(j)}(z)}{z} \right]}.
		\end{aligned}
		\end{equation}
		As $W_{jj} \geq 0$, $\bar W_{(j)} \leq \bar W$, so we find the same equation as Equation \eqref{sqrt_ineq}:
		\begin{equation}\label{}
		\begin{aligned}
			\left| \Re\left[- \frac{\tilde m_{(j)}(z)}{z} \right] \right|& \leq \frac{|u|}{v} \Im\left[- \frac{\tilde m_{(j)}(z)}{z} \right] + \frac{|z|\sqrt{\bar W}}{v\sqrt{v}}  \sqrt{ \Im\left[- \frac{\tilde m_{(j)}(z)}{z} \right]}.
		\end{aligned}
		\end{equation}
		So, from the previous proof of the first inequality, we deduce immediatly the second one, which complete the proof of this lemma:
		\begin{equation}\label{}
		\begin{aligned}
			&u \neq 0 \implies \left \lVert \left( - \frac{\tilde m_{(j)}(z)}{z} A + I \right)^{-1} \right \rVert \leq 16 \left( \frac{|z|^2}{4v^2|u|}\bar W \lVert A \rVert +1 \right) \times \max \left(\frac{1}{3}, \frac{|u|}{v} \right), \\
			&u = 0 \implies \left \lVert \left( - \frac{\tilde m_{(j)}(z)}{z} A + I \right)^{-1} \right \rVert \leq \max\left(2, \frac{4\bar W\lVert A \rVert}{v} \right).
		\end{aligned}
		\end{equation}		
\end{proof}

\begin{crl}\label{tech1_crl1}
	Let $j \in \llbracket 1, N \rrbracket$ and $z = u + iv, v > 0$. Then, for any matrices $A$ and $B$ of same size $n \times n$, $A$ Hermitian non-negative, we have:
	\begin{equation}\label{}
	\begin{aligned}
		&\left | \tr\left[B\left(\left( - \frac{\tilde m_{n}(z)}{z} A + I \right)^{-1} - \left( - \frac{\tilde m_{(j)}(z)}{z} A + I \right)^{-1}\right) \right]\right | \leq 
		\left|\frac{\tilde m_{n}(z)}{z} - \frac{\tilde m_{(j)}(z)}{z} \right|n \lVert B \rVert \lVert A \rVert f(z,\lVert A \rVert)^2.
	\end{aligned}
	\end{equation}
\end{crl}

\begin{proof}
	For any invertible matrices $C_1, C_2$ of the same size than $B$, we have:
	\begin{equation}\label{}
	\begin{aligned}
		|\tr\left[B(C_1^{-1} - C_2^{-1}) \right]| &= \tr\left[BC_1^{-1}(C_2-C_1)C_2^{-1} \right] \\
		|\tr\left[B(C_1^{-1} - C_2^{-1}) \right]| & \leq \lVert B \rVert \times \lVert C_1^{-1} \rVert \times \lVert C_2^{-1} \rVert \times n \lVert C_2-C_1\rVert.
	\end{aligned}
	\end{equation}
	From that point, the result comes immediately from Lemma \ref{tech1}.
\end{proof}

\begin{crl}\label{tech1_crl2}
	Let $j \in \llbracket 1, N \rrbracket$ and $z = u + iv, v > 0$. Then, for any $n \times n$ matrix $A$ Hermitian non-negative, and $r \in \C^n$, $\lVert r \rVert$ denoting its euclidean norm, we have:
	\begin{equation}\label{}
	\begin{aligned}
		&\left | r^* \left(\left( - \frac{\tilde m_{n}(z)}{z} A + I \right)^{-1} - \left( - \frac{\tilde m_{(j)}(z)}{z} A + I \right)^{-1}\right) r\right | \leq 
		\left|\frac{\tilde m_{n}(z)}{z} - \frac{\tilde m_{(j)}(z)}{z} \right| \lVert r \rVert^2 \lVert A \rVert f(z,\lVert A \rVert)^2.
	\end{aligned}
	\end{equation}
\end{crl}

\begin{proof}
	For any invertible matrices $C_1, C_2$ of size $n \times n$, we have:
	\begin{equation}\label{}
	\begin{aligned}
		r^*(C_1^{-1} - C_2^{-1}) r &= r^*C_1^{-1}(C_2-C_1)C_2^{-1} r \\
		r^*(C_1^{-1} - C_2^{-1}) r & \leq \lVert r \rVert^2 \times \lVert C_1^{-1} \rVert \times \lVert C_2^{-1} \rVert \times \lVert C_2-C_1\rVert.
	\end{aligned}
	\end{equation}
	From that point, the result comes immediately from Lemma \ref{tech1}.
\end{proof}

\begin{lem}\label{tech2}
	Let $j \in \llbracket 1, N \rrbracket$. We denote: $A_{(j)} = \sum_{i \neq j} W_{ii} r_i r_i^*.$ Then,
	\begin{equation}\label{}
	\begin{aligned}
		\left|\frac{\tilde m_{n}(z)}{z} - \frac{\tilde m_{(j)}(z)}{z} \right| \leq \frac{2 \log n}{Nv} + \frac{\lVert A_{(j)} \rVert}{|z|vN}.
	\end{aligned}
	\end{equation}
\end{lem}

\begin{proof}
	Let $j \in \llbracket 1, N \rrbracket$. We have, using Lemma \ref{tech0}:
	\begin{equation}\label{}
	\begin{aligned}
		\left|\frac{\tilde m_{n}(z)}{z} - \frac{\tilde m_{(j)}(z)}{z} \right| & = &&\left|\frac{W_{jj}r_j^*(B_n - zI)^{-1}r_j}{zN} + \sum_{i \neq j} \frac{W_{ii}}{zN}  r_i^*\left[(B_n - zI)^{-1} - (B_{(j)} - zI)^{-1} \right]r_i \right| \\
		\left|\frac{\tilde m_{n}(z)}{z} - \frac{\tilde m_{(j)}(z)}{z} \right|& \leq &&\left|\frac{W_{jj}}{zN}\left(1 - \frac{1}{1 + r_j^*(B_{(j)} - zI)^{-1}r_j} \right) \right| \\
		& &&+ \left|\frac{1}{zN} \frac{r_j^*(B_{(j)} - zI)^{-1} A_{(j)} (B_{(j)} - zI)^{-1}r_j}{1+ r_j^* (B_{(j)} - zI)^{-1}r_j} \right|.
	\end{aligned}
	\end{equation}
	For the first term, using Lemma \ref{ineq}, we have:
	\begin{equation}\label{}
	\begin{aligned}
		\left|\frac{W_{jj}}{zN}\left(1 - \frac{1}{1 + r_j^*(B_{(j)} - zI)^{-1}r_j} \right) \right| \leq \frac{\log n}{N|z|} + \frac{\log n}{Nv} \leq \frac{2\log n}{Nv}.
	\end{aligned}
	\end{equation}
	For the second term, we have:
	\begin{equation}\label{}
	\begin{aligned}
		&\left|\frac{1}{zN} \frac{r_j^*(B_{(j)} - zI)^{-1} A_{(j)} (B_{(j)} - zI)^{-1}r_j}{1+ r_j^* (B_{(j)} - zI)^{-1}r_j} \right| \leq  \frac{\lVert A_{(j)}\rVert}{|z|N} \frac{\lVert (B_{(j)} - zI)^{-1}r_j \rVert^2}{|1+ r_j^* (B_{(j)} - zI)^{-1}r_j| } \\
		&\left|\frac{1}{zN} \frac{r_j^*(B_{(j)} - zI)^{-1} A_{(j)} (B_{(j)} - zI)^{-1}r_j}{1+ r_j^* (B_{(j)} - zI)^{-1}r_j} \right| \leq  \frac{\lVert A_{(j)}\rVert}{|z|vN}.
	\end{aligned}
	\end{equation}
	Using the proof of Lemma 2.6 \cite{Silverstein1995b}, we have additionally that:
	\begin{equation}\label{}
	\begin{aligned}
		&\frac{\lVert (B_{(j)} - zI)^{-1}r_j \rVert^2}{|1+ r_j^* (B_{(j)} - zI)^{-1}r_j| }  \leq \frac{1}{v}.
	\end{aligned}
	\end{equation}
	So,
	\begin{equation}\label{}
	\begin{aligned}
		&\left|\frac{1}{zN} \frac{r_j^*(B_{(j)} - zI)^{-1} A_{(j)} (B_{(j)} - zI)^{-1}r_j}{1+ r_j^* (B_{(j)} - zI)^{-1}r_j} \right| \leq  \frac{\lVert A_{(j)}\rVert}{|z|vN},
	\end{aligned}
	\end{equation}
	which concludes the proof.
\end{proof}

\begin{lem}\label{tech3}
	We have:
	\begin{equation}\label{}
	\begin{aligned}
		&\max_{j \leq N}\lVert q_j \rVert^2 \longrightarrow 1 \text{ a.s.}, \\
		&\forall p \in \N, \max_{j \leq N} \frac{\log(n)^p}{N} \lVert A_{(j)} \rVert \longrightarrow 0 \text{ a.s.}
	\end{aligned}
	\end{equation}
\end{lem}

\begin{proof}
	The first convergence is a direct use of Lemma \ref{lemma_rmt} as remarked in \cite{Silverstein1995b} p.338.
	
	Let $p \in \N$. The second convergence comes from the following derivations for $j \in \llbracket 1,N \rrbracket$:
	\begin{equation}\label{}
	\begin{aligned}
		\frac{\log(n)^{2p}}{N^2} \lVert A_{(j)} \rVert^2 &\leq \frac{1}{N^2} \tr \left[\left(\sum_{i \neq j} W_{ii} r_i r_i^* \right) \right] \\
		&\leq \frac{c_n^2\log(n)^{2p}}{N^2} \left(\sum_{i \neq j} W_{ii}^4 \lVert T \rVert^2 \lVert q_i \rVert^4 + \sum_{i \neq j}\sum_{i' \neq i, i' \neq j} W_{ii}^2 W_{i'i'}^2| q_{i'}^* T q_i|^2  \right) \\
		&\leq \frac{c_n^2\log(n)^{2p}}{N^2} \left(\sum_{i = 1}^N W_{ii}^4 \lVert T \rVert^2 \lVert q_i \rVert^4 + \sum_{i = 1}^N \sum_{i' \neq i} W_{ii}^2 W_{i'i'}^2| q_{i'}^* T q_i|^2  \right) \\
		&\leq \frac{c_n^2 \log(n)^{2p+4}}{N^2} \left(\sum_{i = 1}^N \log(n)^2 \lVert q_i \rVert^4 + \sum_{i = 1}^N \sum_{i' \neq i} | q_{i'}^* T q_i|^2  \right) \\
		\frac{1}{N^2} \lVert A_{(j)} \rVert^2 &\leq \frac{c_n^2 \log(n)^{2p+6}}{N} \max_{i \leq N} \lVert q_i \rVert^4 + \frac{c_n^2 \log(n)^{2p+4}}{N} \max_{i \leq N} \sum_{i' \neq i} | q_{i'}^* T q_i|^2.
	\end{aligned}
	\end{equation}
	The upper bound does not depend on $j$ anymore, so we have:
	\begin{equation}\label{}
	\begin{aligned}
		\max_{j \leq N} \frac{1}{N^2} \lVert A_{(j)} \rVert^2 &\leq \frac{c_n^2 \log(n)^{2p+6}}{N} \max_{i \leq N} \lVert q_i \rVert^4 + \frac{c_n^2 \log(n)^{2p+4}}{N} \max_{i \leq N} \sum_{i' \neq i} | q_{i'}^* T q_i|^2.
	\end{aligned}
	\end{equation}
	From the first part of the proof, we have that $ \max_{i \leq N} \lVert q_i \rVert^4 \longrightarrow 1 \text{ a.s.}$, so:
	\begin{equation}\label{}
	\begin{aligned}
		\frac{c_n^2 \log(n)^{2p+6}}{N} \max_{i \leq N} \lVert q_i \rVert^4 \longrightarrow 0 \text{ a.s.}
	\end{aligned}
	\end{equation}
	For the second term, we have:
	\begin{equation}\label{}
	\begin{aligned}
		\frac{c_n^2 \log(n)^{2p+4}}{N} \max_{i \leq N} \sum_{i' \neq i} | q_{i'}^* T q_i|^2 &= \frac{c_n^2 \log(n)^{2p+4}}{N} \max_{i \leq N} \sum_{i' \neq i} q_{i'}^* (T q_i q_i^* T) q_{i'}.
	\end{aligned}
	\end{equation}
	Let $i \in \llbracket 1,N \rrbracket$. We use Lemma \ref{lemma_rmt} in dimension $n(N-1)$ with $Y = (Z_{ki'})_{k,i'\neq i} \in \R^{n \times (N-1)}$ (in vectorized form) and 
	$C = \begin{pmatrix}
		Tq_i q_i^*T & & (0) \\
		& \ddots & \\
		(0) & & Tq_i q_i^*T
	\end{pmatrix}$ of size $n(N-1) \times n(N-1)$. We have then:
	\begin{equation}\label{}
	\begin{aligned}
		&\E\left[\left|\frac{ \log(n)^{2p+4}}{N} \sum_{i' \neq i} q_{i'}^* (T q_i q_i^* T) q_{i'} - \frac{(N-1) \log(n)^{2p+4}}{nN} \lVert T q_i \rVert^2 \right|^6\right] \leq \\
		& \frac{Kn^3(N-2)^3 \log\left(n(N-2)\right)^6 \log(n)^{6(2p+4)}}{N^6n^6} \lVert C \rVert^6 \leq \\
		&\frac{K \log\left(nN\right)^6}{N^3n^3}\log(n)^{6(2p+8)}.
	\end{aligned}
	\end{equation}
	So, for all $\epsilon > 0$:
	\begin{equation}\label{}
	\begin{aligned}
		&\mathbb{P}\left[\left|\max_{i \leq N} \frac{1}{N} \sum_{i' \neq i} q_{i'}^* (T q_i q_i^* T) q_{i'} - \max_{i \leq N}\frac{N-1}{nN} \lVert T q_i \rVert^2 \right| \geq \epsilon\right] \leq \\ 
		&N\times \mathbb{P}\left[\left|\frac{1}{N} \sum_{i' \neq 1} q_{i'}^* (T q_1 q_1^* T) q_{i'} - \frac{N-1}{nN} \lVert T q_1 \rVert^2 \right| \geq \epsilon\right] \leq \\
		&\leq \frac{K \log\left(nN\right)^6}{\epsilon^6N^2n^3} \log(n)^{6(2p+8)},
	\end{aligned}
	\end{equation}
	which is summable. So,
	\begin{equation}\label{}
	\begin{aligned}
		\max_{i \leq N} \frac{\log(n)^{2p+4}}{N} \sum_{i' \neq i} q_{i'}^* (T q_i q_i^* T) q_{i'} - \max_{i \leq N}\frac{(N-1)\log(n)^{2p+4}}{nN} \lVert T q_i \rVert^2 \longrightarrow 0 \text{ a.s.}
	\end{aligned}
	\end{equation}
	And, from the first part of the proof:
	\begin{equation}\label{}
	\begin{aligned}
		\left| \max_{i \leq N}\frac{c_n^2(N-1)\log(n)^{2p+4}}{nN} \lVert T q_i \rVert^2 \right| \leq \frac{c_n^2 \log(n)^{2p+6}}{nN} \max_{i \leq N} \lVert q_i \rVert^2 \longrightarrow 0 \text{ a.s.}
	\end{aligned}
	\end{equation}
	We can conclude the proof:
	\begin{equation}\label{}
	\begin{aligned}
		&\max_{j \leq N} \frac{\log(n)^p}{N} \lVert A_{(j)} \rVert \longrightarrow 0 \text{ a.s.}
	\end{aligned}
	\end{equation}
\end{proof}

\subsubsection{Proof that $\max_{j \leq N} |d_j| \longrightarrow 0$ a.s.} 
\paragraph{Proof that $\max_{j \leq N} |d_j^{(1)}| \longrightarrow 0$ a.s.}
Let $j \in \llbracket 1, N \rrbracket$. We recall that:
\begin{equation}\label{}
\begin{aligned}
	d_j^{(1)} &= W_{jj}q_j^* T_n^{1/2} (B_{(j)} - zI)^{-1} \left[\left(-\frac{\tilde{m}_n(z)}{z}  T_n + I\right)^{-1} - \left(-\frac{\tilde{m}_{(j)}(z)}{z}  T_n + I\right)^{-1} \right]T_n^{1/2}q_j.
\end{aligned}
\end{equation}
Then, using Corollary \ref{tech1_crl2} and Lemma \ref{ineq}, we have:
\begin{equation}\label{}
\begin{aligned}
	|d_j^{(1)}| &\leq W_{jj} \lVert (B_{(j)} - zI)^{-1} \rVert \lVert T \rVert^2 \lVert q_j \rVert^2 \left|\frac{\tilde m_{n}(z)}{z} - \frac{\tilde m_{(j)}(z)}{z} \right| f(z,\log(n))^2 \\
	|d_j^{(1)}| & \leq \frac{\log(n)^3}{v}f(z,\log(n))^2 \lVert q_j \rVert^2 \left|\frac{\tilde m_{n}(z)}{z} - \frac{\tilde m_{(j)}(z)}{z} \right|.
\end{aligned}
\end{equation}
Using Lemma \ref{tech2}, we have:
\begin{equation}\label{}
\begin{aligned}
	|d_j^{(1)}| & \leq \frac{\log(n)^3}{v} f(z,\log(n))^2  \left(\frac{2 \log n}{Nv} + \frac{\lVert A_{(j)} \rVert}{|z|vN}\right)\lVert q_j \rVert^2.
\end{aligned}
\end{equation}
By Assumption (d), $\int x dF^{W_n}(x) \longrightarrow \int x dD(x) < \infty$ a.s., \textit{i.e.} $\bar W \longrightarrow \int x dD(x) \leq \infty$ a.s., so $\bar W$ is bounded a.s. As a consequence, $f(z,\log(n))^2 = O(\log(n)^2)$ a.s.. Finally, using Lemma \ref{tech3}, we can conclude:
\begin{equation}\label{}
\begin{aligned}
	\max_{j \leq N} |d_j^{(1)}| \longrightarrow 0 \text{ a.s.}
\end{aligned}
\end{equation}

\paragraph{Proof that $\max_{j \leq N} |d_j^{(2)}| \longrightarrow 0$ a.s.}
Let $j \in \llbracket 1, N \rrbracket$. We recall that:
\begin{equation}\label{}
\begin{aligned}
	d_j^{(2)} = & W_{jj}q_j^* T_n^{1/2} (B_{(j)} - zI)^{-1}\left(-\frac{\tilde{m}_{(j)}(z)}{z}  T_n + I\right)^{-1}T_n^{1/2}q_j \\
	& - \frac{W_{jj}}{n}\tr \left( \left(-\frac{ \tilde{m}_{(j)}(z)}{z} T_n + I\right)^{-1} T_n (B_{(j)} - zI)^{-1}\right).
\end{aligned}
\end{equation}
Using Lemma \ref{lemma_rmt}, we have:
\begin{equation}\label{}
\begin{aligned}
	\E\left[\left|d_j^{(2)}\right|^6\right] & \leq \frac{K \log(n)^{12}}{n^3} \left \lVert W_{jj} \left(-\frac{ \tilde{m}_{(j)}(z)}{z} T_n + I\right)^{-1} T_n (B_{(j)} - zI)^{-1} \right \rVert^6.
\end{aligned}
\end{equation}
Using Lemmas \ref{ineq} and \ref{tech1}, we have:
\begin{equation}\label{}
\begin{aligned}
	\E\left[\left|d_j^{(2)}\right|^6\right] & \leq \frac{K \log(n)^{24}}{v^6n^3} f(z, \log(n))^6.
\end{aligned}
\end{equation}
As argued above, $f(z, \log(n)) = O(\log(n))$ a.s., so $\sum_{j=1}^N \E\left[\left|d_j^{(2)}\right|^6\right]$ is a.s. summable. So, we can conclude:
\begin{equation}\label{}
\begin{aligned}
	\max_{j \leq N} |d_j^{(2)}| \longrightarrow 0 \text{ a.s.}
\end{aligned}
\end{equation}

\paragraph{Proof that $\max_{j \leq N} |d_j^{(3)}| \longrightarrow 0$ a.s.}
Let $j \in \llbracket 1, N \rrbracket$. We recall that:
\begin{equation}\label{}
\begin{aligned}
	d_j^{(3)} &= \frac{W_{jj}}{n}\tr \left(  \left[\left(-\frac{\tilde{m}_{(j)}(z)}{z}  T_n + I\right)^{-1} - \left(-\frac{\tilde{m}_{n}(z)}{z}  T_n + I\right)^{-1} \right] T_n (B_{(j)} - zI)^{-1}\right).
\end{aligned}
\end{equation}
Using Corollary \ref{tech1_crl1} and Lemma \ref{ineq}, we have:
\begin{equation}\label{}
\begin{aligned}
	|d_j^{(3)}| &\leq \frac{\log(n)^3}{nv} \left|\frac{\tilde m_{n}(z)}{z} - \frac{\tilde m_{(j)}(z)}{z} \right|n f(z, \log(n))^2.
\end{aligned}
\end{equation}
Using Lemma \ref{tech2}, we have:
\begin{equation}\label{}
\begin{aligned}
	|d_j^{(3)}| &\leq \frac{\log(n)^3}{v} \left(\frac{2 \log n}{nv} + \frac{\lVert A_{(j)} \rVert}{|z| v n}\right) f(z, \log(n))^2.
\end{aligned}
\end{equation}
As we argued previously, $f(z, \log(n))^2 = O(\log(n)^2)$ a.s., so from Lemma \ref{tech3} we can conclude:
\begin{equation}\label{}
\begin{aligned}
	\max_{j\ \leq N} |d_j^{(3)}|\longrightarrow 0 \text{ a.s.}
\end{aligned}
\end{equation}

\paragraph{Proof that $\max_{j \leq N} |d_j^{(4)}| \longrightarrow 0$ a.s.}
Let $j \in \llbracket 1, N \rrbracket$. We recall that:
\begin{equation}\label{}
\begin{aligned}
	d_j^{(4)} &= \frac{W_{jj}}{n}\tr \left(\left(-\frac{\tilde{m}_{n}(z)}{z}  T_n + I\right)^{-1} T_n \left[(B_{(j)} - zI)^{-1} - (B_{n} - zI)^{-1}\right]\right) \\
\end{aligned}
\end{equation}
Using Lemma \ref{tech0}, we have:
\begin{equation}\label{}
\begin{aligned}
	d_j^{(4)} &= \frac{W_{jj}}{n} \frac{r_j^* (B_{(j)} - zI)^{-1} \left(-\frac{\tilde{m}_{n}(z)}{z}  T_n + I\right)^{-1} T  (B_{(j)} - zI)^{-1}r_j}{1 + r_j^* (B_{(j)} - zI)^{-1} r_j}
\end{aligned}
\end{equation}
So,
\begin{equation}\label{}
\begin{aligned}
	|d_j^{(4)}| &= \frac{\log(n)}{n} \left \lVert \left(-\frac{\tilde{m}_{n}(z)}{z}  T_n + I\right)^{-1} T \right \rVert \frac{\lVert (B_{(j)} - zI)^{-1}r_j \rVert^2}{|1 + r_j^* (B_{(j)} - zI)^{-1} r_j|}
\end{aligned}
\end{equation}
Using the proof of Lemma 2.6 \cite{Silverstein1995b}, we have that:
\begin{equation}\label{}
\begin{aligned}
	&\frac{\lVert (B_{(j)} - zI)^{-1}r_j \rVert^2}{|1+ r_j^* (B_{(j)} - zI)^{-1}r_j| }  \leq \frac{1}{v}.
\end{aligned}
\end{equation}
So, using Lemma \ref{tech2}, we have:
\begin{equation}\label{}
\begin{aligned}
	|d_j^{(4)}| &= \frac{\log(n)^2 f(z, \log(n))}{vn}
\end{aligned}
\end{equation}
As argued before, $f(z, \log(n)) = O(\log(n))$ a.s., so we can conclude:
\begin{equation}\label{}
\begin{aligned}
	\max_{j \leq N} |d_j^{(4)}| \longrightarrow 0 \text{ a.s.}
\end{aligned}
\end{equation}

We can now conclude this section. The last four points prove  that:
\begin{equation}\label{}
\begin{aligned}
	\max_{j \leq N} |d_j| \longrightarrow 0 \text{ a.s.}
\end{aligned}
\end{equation}
And from Lemma \ref{ineq}, we have for $j \in \llbracket 1, N \rrbracket$:
\begin{equation}\label{}
\begin{aligned}
	&\frac{1}{|z(1 + r_j^* (B_{(j)} -zI)^{-1}r_j)|} \leq \frac{1}{v}.
\end{aligned}
\end{equation}
So, 
\begin{equation}\label{}
\begin{aligned}
	&\frac{1}{N} \sum_{j=1}^N \frac{-1}{z(1 + r_j^* (B_{(j)} -zI)^{-1}r_j)} d_j \underset{n \rightarrow \infty}{\longrightarrow} 0 \text{ a.s.}
\end{aligned}
\end{equation}
Using Equation \eqref{eq_cvg}, we can now conlude that:
\begin{equation}\label{cvg}
\begin{aligned}
	&\frac{1}{n} \tr\left((\tilde{m}_n(z) T_n - zI)^{-1}\right) - m_n(z) \underset{n \rightarrow \infty}{\longrightarrow} 0 \text{ a.s.}
\end{aligned}
\end{equation} 

\subsection{Convergences and functional equation} 
For this section, we introduced an object used in \cite{Ledoit2009}, we define:
\begin{equation}\label{}
\begin{aligned}
	\Theta_n^{(1)}(z) =\frac{1}{n} \tr \left((\mathrm{B}_n - zI)^{-1}T_n \right).
\end{aligned}
\end{equation}

\subsubsection{Proof that a.s., $\forall p \in \N, \log(n)^p \left|\tilde m_n(z) \Theta_{n}^{(1)}(z)- (1 + zm_n(z)) \right| \underset{n \rightarrow \infty}{\longrightarrow} 0$} 
Let $p \in \N$. We have:
\begin{equation}\label{}
\begin{aligned}
	\left|\tilde m_{n}(z) \Theta_{n}^{(1)}(z) - \left(1 + z m_{n}(z)\right)\right| &= \left|\tilde m_{n}(z) \Theta_{n}^{(1)}(z) - \frac{1}{n} \sum_{j=1}^N \frac{r_j^*(B_{(j)}-zI)^{-1}r_j}{1+r_j^*(B_{(j)}-zI)^{-1}r_j}\right| \\
	&= \left|\frac{1}{N} \sum_{j=1}^N \frac{W_{jj} \left(\Theta_{n}^{(1)}(z) - q_j^*(B_{(j)}-zI)^{-1}q_j\right)}{1+r_j^*(B_{(j)}-zI)^{-1}r_j}\right| \\
	\left|\tilde m_{n}(z) \Theta_{n}^{(1)}(z) - \left(1 + z m_{n}(z)\right)\right| & \leq \frac{|z|}{v}\max_{j \leq N}  \log(n) \left|\Theta_{n}^{(1)}(z) - q_j^*T^{1/2}(B_{(j)}-zI)^{-1}T^{1/2}q_j \right|.
\end{aligned}
\end{equation} 
Using Lemma \ref{lemma_rmt}, we have that:
\begin{equation}\label{}
\begin{aligned}
	\E\left[\left|\log(n)^{p}\Theta_{n}^{(1)}(z) - \log(n)^{p}q_j^*T^{1/2}(B_{(j)}-zI)^{-1}T^{1/2}q_j \right|^6\right] \leq K \frac{\log(n)^{6(p+3)}}{n^3 v^6}.
\end{aligned}
\end{equation} 
So, 
\begin{equation}\label{}
\begin{aligned}
	\max_{j \leq N} \left|\log(n)^{p}\Theta_{n}^{(1)}(z) - \log(n)^{p}q_j^*T^{1/2}(B_{(j)}-zI)^{-1}T^{1/2}q_j \right| \underset{n \rightarrow \infty}{\longrightarrow} 0 \text{ a.s.}
\end{aligned}
\end{equation} 
As $\N$ is countable, we have also:
\begin{equation}\label{theta_cvg}
\begin{aligned}
	\text{a.s., } \forall p \in \N, \max_{j \leq N} \log(n)^{p}\left|\Theta_{n}^{(1)}(z) - q_j^*T^{1/2}(B_{(j)}-zI)^{-1}T^{1/2}q_j \right| \underset{n \rightarrow \infty}{\longrightarrow} 0 \text{ a.s.}
\end{aligned}
\end{equation} 
We can conclude:
\begin{equation}\label{theta_cvg2}
\begin{aligned}
	\text{a.s., } \forall p \in \N, \log(n)^p \left|\tilde m_n(z) \Theta_{n}^{(1)}(z)- (1 + zm_n(z)) \right| \underset{n \rightarrow \infty}{\longrightarrow} 0.
\end{aligned}
\end{equation} 

From Assumption (d), $\bar W \rightarrow \int xdD(x) \in \R_+$ a.s. We focus now on trajectories where $\bar W \rightarrow \int xdD(x) \in \R_+$, Equations \eqref{cvg} and \eqref{theta_cvg} hold, $F^{T_n} \implies H$ and $F^{W_n} \implies D$.

Then, $(\tilde m_n(z))$ is bounded. Indeed, from Lemma \ref{ineq}, $| \tilde m_n(z) | \leq \bar W \frac{|z|}{v}$. Then, as $\Im[\tilde m_n(z)] \leq 0$,it exists a subsequence $\{n_i\}$ of $\N$ and $\tilde m(z) \in \C \backslash \C_+$ such that $\tilde m_{n_i}(z) \underset{i \rightarrow \infty}{\longrightarrow} \tilde m(z) \in \C \backslash \C_+$. 

\subsubsection{Proof that $m_{n_i}(z) \underset{i \rightarrow \infty}{\longrightarrow} \int \frac{1}{\tau \tilde m(z) - z}dH(\tau)$.} 
We want to prove that:
\begin{equation}\label{}
\begin{aligned}
	&m_{n_i}(z) -  \int \frac{1}{\tau \tilde m(z) - z}dH(\tau) \underset{i \rightarrow \infty}{\longrightarrow} 0.
\end{aligned}
\end{equation} 
From Equation \eqref{cvg}, it is equivalent to prove that:
\begin{equation}\label{}
\begin{aligned}
	&\int \frac{1}{\tau \tilde m_{n_i}(z) - z}dF^{T_{n_i}}(\tau) -  \int \frac{1}{\tau \tilde m(z) - z}dH(\tau) \underset{i \rightarrow \infty}{\longrightarrow} 0.
\end{aligned}
\end{equation} 

We prove that $\int \frac{1}{\tau \tilde m_{n_i}(z) - z}dF^{T_{n_i}}(\tau) \underset{i \rightarrow \infty}{\longrightarrow} -\frac{1}{z}$ using the Lebesgue's convergence theorem for weakly converging measures, as detailed in Corollary 5.1 \cite{Feinberg2019}. We denote: $f: \tau \in \R \rightarrow \frac{1}{\tau \tilde m(z) - z}$ and $f_{i}: \tau \in \R \rightarrow \frac{1}{\tau \tilde m_{n_i}(z) - z}$. Regarding the hypotheses of the theorem, we have:
\begin{itemize}
	\item $(f_i)_i$ is a.u.i. w.r.t. $(F^{T_{n_i}})_i$ (see (2.4) \cite{Feiberg2019} for a definition). Indeed, $\forall \tau \in \R_+, \forall i \in \R_+, |f_i(\tau)| \leq \frac{1}{v}$, so $\lim_{K \rightarrow +\infty} \limsup_{i \rightarrow +\infty} \int |f_i(\tau)| 1_{[K,+\infty[}(|f_i(\tau)|) dF^{T_{n_i}} = 0$.
	\item $F^{T_{n_i}} \implies H$ by assumption.
	\item Let $\tau \in \R_+$ and $\epsilon > 0$. By assumption, $\bar W_{n_i}$ is bounded and we denote by $\kappa$ one of its finite upper bound. By \eqref{ineq}, we have that $\forall i, |\tilde m_{n_i}(z) | \leq \frac{\kappa}{|z|v}$, so $\tilde m(z) \leq \frac{\kappa}{|z|v}$. We define:
		\begin{equation}\label{}
		\begin{aligned}
			\delta := \min\left(1,\frac{v^2}{\frac{\kappa}{|z|v} + \tau + 1}\right) > 0.
		\end{aligned}
		\end{equation} 
	Then, there exists $i_0 \in \N$ such that $\forall i \geq i_0,  |\tilde m_{n_i}(z) - \tilde m(z) | \leq \delta$. Now, let $\tau' \in ]\tau - \delta, \tau + \delta[ \cap \R_+$ and $i \geq i_0$. Then,
	\begin{equation}\label{}
	\begin{aligned}
		\left|f_i(\tau') - f(\tau) \right| &= \frac{|\tau' \tilde m_{n_i}(z) - \tau \tilde m(z)|}{|\tau' \tilde m_{n_i}(z) - z|\times|\tau \tilde m(z) - z|} \\
		&\leq \frac{1}{v^2}\left(|(\tau' - \tau) \tilde m(z)| + |\tau' (\tilde m_{n_i}(z)-\tilde m(z))| \right) \\
		&\leq \frac{\delta}{v^2}\left(\frac{\kappa}{|z|v} + \tau + 1 \right) \\
		\left|f_i(\tau') - f(\tau) \right| &\leq \epsilon.
	\end{aligned}
	\end{equation} 
	So $\lim_{i \rightarrow \infty, \tau' \rightarrow \tau} f_i(\tau')$ exists.
\end{itemize}
So, using Corollary 5.1 \cite{Feinberg2019} on the real and imaginary parts, we deduce that $\lim_{i \rightarrow \infty} \int f_i(\tau) dF^{T_{n_i}}(\tau)$ exists and:
\begin{equation}\label{tau}
\begin{aligned}
	&\lim_{i \rightarrow \infty} \int f_i(\tau) dF^{T_{n_i}}(\tau) = \int f(\tau) dH(\tau).
\end{aligned}
\end{equation} 

It immediatly leads to:
\begin{equation}\label{tau}
\begin{aligned}
	&m_{n_i}(z) -  \int \frac{1}{\tau \tilde m(z) - z}dH(\tau) \underset{i \rightarrow \infty}{\longrightarrow} 0.
\end{aligned}
\end{equation} 

\subsubsection{Case $H = 1_{[0,+\infty[}$.}
Suppose $H = 1_{[0,+\infty[}$. We have then that a.s. $m_n(z) \underset{n \rightarrow +\infty}{\longrightarrow} -\frac{1}{z}$. Moreover, the equation $m = \int \frac{\delta}{1 + \delta c \int \frac{\tau}{\tau m - z} dH(\tau)}dD(\delta) = \int \delta dD(\delta)$ has trivially a unique solution in $\C \backslash \C_+$ that we denote $\tilde m^{(0)}(z) :=  \int \delta dD(\delta)$. As $-\frac{1}{z} = \int \frac{1}{\tau \tilde m^{(0)}(z) - z}dH(\tau)$, we indeed have that $m_n(z) \underset{n \rightarrow +\infty}{\longrightarrow} \int \frac{1}{\tau \tilde m^{(0)}(z) - z}dH(\tau)$. And $z \mapsto -\frac{1}{z}$ is the Cauchy-Stieltjes transform of the p.d.f $1_{[0,\infty[}$ which complete the proof of the Theorem in the case $H = 1_{[0,+\infty[}$.

Until the end of the proof, we now suppose that $H(]0,+\infty[) > 0$.

\subsubsection{Case $D = 1_{[0,+\infty[}$.}
Suppose $D = 1_{[0,+\infty[}$. Then, $\bar W_n \underset{n \rightarrow +\infty}{\longrightarrow} 0$. By \eqref{ineq}, $\tilde m_n(z) \leq \frac{\bar W}{|z|v}$. So the complete sequence (not the subsequence) $\tilde m_n(z) \underset{n \rightarrow +\infty}{\longrightarrow} 0$. From the previous section on Equation \eqref{tau}, we proved that if $\tilde m_n(z)$ converges to some $\tilde m(z)$, then $m_n(z) -  \int \frac{1}{\tau \tilde m(z) - z}dH(\tau) \underset{n \rightarrow \infty}{\longrightarrow} 0.$ So, with $\tilde m^{(1)}(z) := 0$, we have that $m_n(z) \underset{n \rightarrow \infty}{\longrightarrow} \int \frac{1}{\tau \tilde m^{(1)}(z) - z}dH(\tau) = - \frac{1}{z}$. Moreover, we remark that $\tilde m^{(1)}(z) = 0$ is the unique solution to the equation $m = \int \frac{\delta}{1 + \delta c \int \frac{\tau}{\tau m - z} dH(\tau)}dD(\delta) = \int \delta dD(\delta)$ in $\C \backslash \C_+$. And $z \mapsto -\frac{1}{z}$ is the Cauchy-Stieltjes transform of the p.d.f $1_{[0,\infty[}$ which complete the proof of the Theorem in the case $D = 1_{[0,+\infty[}$.

Until the end of the proof, we now suppose that $D(]0,+\infty[) > 0$.

\subsubsection{Case $\tilde m(z) = 0$.}
Suppose $\tilde m(z) = 0$. Then, $m_{n_i}(z) \underset{i \rightarrow \infty}{\longrightarrow} - \frac{1}{z}$ and by Cauchy-Stieltjes transform property $F^{B_{n_i}} \implies 1_{[0,+\infty[}$. Consequently, $\E\left[\frac{1}{n_i} \tr(B_{n_i})\right] \underset{i \rightarrow +\infty}{\longrightarrow} 0$. As $\E\left[\frac{1}{n_i} \tr(B_{n_i})\right] = \frac{1}{n_i}\tr(T_{n_i}) \times \bar W_{n_i}$. $\bar W_{n_i}$ converges by assumption to $\int \delta dD(\delta)$, and we supposed $D(]0,+\infty[) > 0$ so $\int \delta dD(\delta) > 0$. We deduce that $\frac{1}{n_i}\tr(T_{n_i}) \underset{i \rightarrow +\infty}{\longrightarrow} 0$. So, $H = 1_{[0,+\infty[}$, which is absurd. So $\tilde m(z) \neq 0$.

\subsubsection{Proof that $1 + zm_{n_i}(z) \underset{i \rightarrow \infty}{\longrightarrow} \int \frac{\delta \Theta^{(1)}(z)}{1 + \delta c \Theta^{(1)}(z)}dD(\delta)$ with $\Theta^{(1)}(z) := \frac{1+zm(z)}{\tilde m(z)}$.} 
We remind that we suppose now: $H(]0,+\infty[) > 0$, $D(]0,+\infty[) > 0$ and $\tilde m_{n_i}(z) \underset{i \rightarrow +\infty}{\longrightarrow} \tilde m(z) \neq 0$.

Firstly, we denote: $M_{n} := \log(n) \max_{j \leq N} \left|\Theta_{n}^{(1)}(z) - \log(n)^{p}q_j^*T^{1/2}(B_{(j)}-zI)^{-1}T^{1/2}q_j \right|$. There is $n_0$ large enough so that $\forall n \geq n_0, M_n \leq \frac{v}{2c_n}$, due to Equation \eqref{theta_cvg}. Then, we have for $n_i \geq n_0$:
\begin{equation}\label{cvg_delta}
\begin{aligned}
	&\left| \left(-c_n(1 + zm_{n_i}(z)) +1\right) - \frac{1}{N}\sum_{j=1}^N \frac{1}{1 + W_{jj} c_{n_i} \Theta_{n_i}^{(1)}} \right| = \\
	&\left| \frac{1}{N}\sum_{j=1}^N \frac{1}{1 + r_j^*(B_{(j)}-zI)^{-1}r_j} - \frac{1}{1 + W_{jj} c_{n_i} \Theta_{n_i}^{(1)}} \right| \leq \\
	& \frac{c_{n_i} \bar W M_{n_i}}{v(v - c_{n_i} M_{n_i})} \underset{i \rightarrow \infty}{\longrightarrow} 0.
\end{aligned}
\end{equation} 
Then, we denote $\Theta^{(1)}(z) := \frac{1+zm(z)}{\tilde m(z)} = \int \frac{\tau}{\tau \tilde m(z) - z}dH(\tau)$. From Equation \eqref{theta_cvg2}, $\Theta^{(1)}_{n_i}\underset{i \rightarrow \infty}{\longrightarrow}  \Theta^{(1)}(z)$. We define:
\begin{equation}\label{}
\begin{aligned}
	&g: \delta \in \R \rightarrow \frac{1}{1 + \delta c \Theta^{(1)}(z)}.
\end{aligned}
\end{equation} 
We want to prove that:
\begin{equation}\label{}
\begin{aligned}
	&\frac{1}{N}\sum_{j=1}^N \frac{1}{1 + W_{jj} c_{n_i} \Theta_{n_i}^{(1)}} \underset{i \rightarrow \infty}{\longrightarrow} \int g(\delta) dD(\delta).
\end{aligned}
\end{equation} 
Let $\delta \in \R$. Remark that, as $\Im[\tilde m(z)]  \leq 0$:
\begin{equation}\label{im_theta}
\begin{aligned}
	&\Im[\Theta^{(1)}(z)] = \Im\left[\frac{1+zm(z)}{\tilde m(z)} \right] = \Im \left[\int \frac{\tau}{\tau \tilde m(z) - z}dH(\tau) \right] \geq v\int \frac{\tau}{|\tau \tilde m(z) - z|^2}dH(\tau) > 0.
\end{aligned}
\end{equation} 
So, from \cite{Silverstein1995c} p.338, we have:
\begin{equation}\label{ineq_delta1}
\begin{aligned}
	&|g(\delta)| \leq \frac{|\Theta^{(1)}(z)|}{\Im[\Theta^{(1)}(z)]}.
\end{aligned}
\end{equation} 
And, for $n_i$ large enough so that $\forall n_i \geq n_1, \Im[\Theta^{(1)}_{n_i}] > 0$:
\begin{equation}\label{ineq_delta2}
\begin{aligned}
	&\left| \frac{1}{1 + \delta c_{n_i} \Theta_{n_i}^{(1)}} - g(\delta)\right| \leq \frac{1}{c_{n_i} \Im[\Theta_{n_i}^{(1)}(z)]}\times \frac{|\Theta^{(1)}(z)|}{\Im[\Theta^{(1)}(z)]} \left|c\Theta^{(1)}(z) - c_{n_i} \Theta^{(1)}_{n_i}\right| .
\end{aligned}
\end{equation} 
Using \eqref{cvg_delta}, \eqref{ineq_delta1} and \eqref{ineq_delta2}, we have that:
\begin{equation}\label{}
\begin{aligned}
	&c_{n_i}\left(1 + zm_{n_i}(z)\right) - 1 \underset{i \rightarrow \infty}{\longrightarrow} \int \frac{ 1}{1 + \delta c \Theta^{(1)}(z)}dD(\delta).
\end{aligned}
\end{equation} 
We conclude that:
\begin{equation}\label{delta}
\begin{aligned}
	&1 + zm_{n_i}(z) \underset{i \rightarrow \infty}{\longrightarrow} \int \frac{ \delta \Theta^{(1)}(z)}{1 + \delta c \Theta^{(1)}(z)}dD(\delta).
\end{aligned}
\end{equation} 
Finally, using \eqref{tau} and \eqref{delta}, we have that:
\begin{equation}\label{tilde_m}
\begin{aligned}
	&\tilde m(z) = \int \frac{\delta}{1 + \delta c \int \frac{\tau}{\tau \tilde m(z) - z} dH(\tau)}dD(\delta).
\end{aligned}
\end{equation} 

\subsection{Unicity}\label{uni} 
In this section, we show that there is at most one solution $m \in \C \backslash \C_+$ so that:
\begin{equation}\label{fun_eq}
\begin{aligned}
	&m = \int \frac{\delta}{1 + \delta c \int \frac{\tau}{\tau m - z} dH(\tau)}dD(\delta).
\end{aligned}
\end{equation}

Let $m \in \C \backslash \C_+$, verifying \eqref{fun_eq}. Then, as $H(]0, +\infty[) > 0$:
\begin{equation}\label{uni_imag}
\begin{aligned}
	\Im[m] &= -\int \frac{\delta^2 c \int \frac{\tau(\tau \Im[-m] + v)}{|\tau m - z|^2}dH(\tau)}{\left| 1 + \delta c \int \frac{\tau}{\tau m - z}dH(\tau) \right|^2} dD(\delta) < 0.
\end{aligned}
\end{equation} 
So, with $\C_- = \{ y \in \C | \Im[y] < 0\}$:
\begin{equation}\label{}
\begin{aligned}
	\left[m \in \C \backslash \C_+ \text{ verifies \eqref{fun_eq} and } H(]0, +\infty[) > 0\right] \iff \left[m \in \C_- \text{ verifies \ref{fun_eq} and } H(]0, +\infty[) > 0\right].
\end{aligned}
\end{equation} 
Let $m_1 = u_1 - i v_1 \in \C_-$ and $m_2 = u_2 - i v_2 \in \C_-$ solving \eqref{fun_eq} at $z = u + iv \in \C_+$. Then, 
\begin{equation}\label{uni_diff}
\begin{aligned}
	m_1 - m_2 &= \int \frac{\delta^2 c \left(\int \frac{\tau}{\tau m_2 - z} dH(\tau) - \int \frac{\tau}{\tau m_1 - z} dH(\tau) \right)}{\left(1 + \delta c \int \frac{\tau}{\tau m_1 - z} dH(\tau)\right)\left(1 + \delta c \int \frac{\tau}{\tau m_2 - z} dH(\tau)\right)} dD(\delta) \\
	m_1 - m_2 &= (m_1 - m_2) \times \int \frac{\delta^2 c \int \frac{\tau^2}{(\tau m_1 - z)(\tau m_2 - z)} dH(\tau) }{\left(1 + \delta c \int \frac{\tau}{\tau m_1 - z} dH(\tau)\right)\left(1 + \delta c \int \frac{\tau}{\tau m_2 - z} dH(\tau)\right)} dD(\delta)
\end{aligned}
\end{equation} 
Using Hölder inequality on the last term and using \eqref{uni_imag} at the end, we have:
\begin{equation}\label{}
\begin{aligned}
	&\left| \int \frac{\delta^2 c \int \frac{\tau^2}{(\tau m_1 - z)(\tau m_2 - z)} dH(\tau) }{\left(1 + \delta c \int \frac{\tau}{\tau m_1 - z} dH(\tau)\right)\left(1 + \delta c \int \frac{\tau}{\tau m_2 - z} dH(\tau)\right)} dD(\delta) \right|^2 \leq \\
	& \left| \int \frac{\delta^2 c \int \frac{\tau^2}{|\tau m_1 - z|^2} dH(\tau) }{\left|1 + \delta c \int \frac{\tau}{\tau m_1 - z} dH(\tau)\right|^2} dD(\delta) \right| \times \left|\int \frac{\delta^2 c \int \frac{\tau^2}{|\tau m_2 - z|^2} dH(\tau) }{\left|1 + \delta c \int \frac{\tau}{\tau m_2 - z} dH(\tau)\right|^2} dD(\delta) \right| < \\
	& \left| \int \frac{\delta^2 c \int \frac{\tau^2 + \tau \frac{v}{v_1}}{|\tau m_1 - z|^2} dH(\tau) }{\left|1 + \delta c \int \frac{\tau}{\tau m_1 - z} dH(\tau)\right|^2} dD(\delta) \right| \times \left|\int \frac{\delta^2 c \int \frac{\tau^2+ \tau \frac{v}{v_2}}{|\tau m_2 - z|^2} dH(\tau) }{\left|1 + \delta c \int \frac{\tau}{\tau m_2 - z} dH(\tau)\right|^2} dD(\delta) \right| = \\
	& \left|\frac{v_1}{v_1} \right| \times \left|\frac{v_2}{v_2} \right| = 1.
\end{aligned}
\end{equation} 
Remark that the inequality is strict because we supposed $H(]0, +\infty[) > 0$. Injecting this inequation in \eqref{uni_diff}, we find:
\begin{equation}\label{}
\begin{aligned}
	|m_1 - m_2| \neq 0 \implies |m_1 - m_2| < |m_1 - m_2|.
\end{aligned}
\end{equation} 
So there is at most one solution $m \in \C \backslash \C_+$ verifying \eqref{fun_eq}.

\subsection{Convergences of $\tilde m_n(z)$, $m_n(z)$ and $\Theta^{(1)}_n(z)$.}
Backing up, we proved that almost surely, $\tilde m_n(z)$ is bounded and every convergent subsequence of $\tilde m_n(z)$ converge towards the unique $\tilde m(z) \in \C \backslash \C_+$ verifying \eqref{tilde_m}. So, a.s., $\tilde m_n(z) \underset{n \rightarrow \infty}{\longrightarrow} \tilde m(z)$.

We also proved that, almost surely, if $\tilde m_n(z) \underset{n \rightarrow \infty}{\longrightarrow} \tilde m(z)$ then:
\begin{itemize}
	\item $m_n(z) \underset{n \rightarrow \infty}{\longrightarrow} m(z) = \int \frac{1}{\tau \tilde m(z) - z} dH(\tau)$,
	\item $\Theta^{(1)}_n(z) \underset{n \rightarrow \infty}{\longrightarrow} \Theta^{(1)}(z) = \int \frac{\tau}{\tau \tilde m(z) - z} dH(\tau)$.
\end{itemize}

So, almost surely, $m_n(z) \underset{n \rightarrow \infty}{\longrightarrow} m(z) = \int \frac{1}{\tau \tilde m(z) - z} dH(\tau)$ where $\tilde m(z)$ is the unique solution in $\C \backslash \C_+$ to \eqref{tilde_m}.

\subsection{$m$ is Stieltjes transform of a p.d.f.}
The last point to prove is that $m$ is Stieltjes transform of a p.d.f. As pointwise limit of Stieltjes transform, it is enough to prove that, using Theorem 7 \cite{Najim2016}, $i y m(iy) \underset{y \rightarrow \infty}{\longrightarrow} -1$ to show that $m$ is a Stieltjes transform of a p.d.f. 

Consider a trajectory where $\bar W$ is bounded, say by $\kappa$, and $\tilde m_n(z) \rightarrow \tilde m(z)$. Then, $\forall z \in \C_+, | \tilde m_n(z)| \leq \kappa \frac{|z|}{v}$. So, $\forall y \in \R_+^*, |\tilde m(iy)| \leq \kappa$. Consequently,
\begin{equation}\label{uni_diff}
\begin{aligned}
	\left | -i y m(iy) - 1 \right| = \left | \int \frac{-\tau \tilde m(iy)}{\tau \tilde m(iy) - iy}dH(\tau) \right| \leq \frac{\kappa}{y}.
\end{aligned}
\end{equation} 
So, $i y m(iy) \underset{y \rightarrow \infty}{\longrightarrow} -1$, which proves that $m$ is a Stieltjes transform a p.d.f. and concludes the proof.

\bibliography{aistats_2025}

\end{document}